\documentclass[final,11pt,oneside]{article}%
\usepackage{amssymb}
\usepackage{amsmath}
\usepackage{amsfonts}
\usepackage{proc2e}
\usepackage{color}
\usepackage{hyperref}%
\usepackage{graphicx}

\definecolor{rltred}{rgb}{0.75,0,0}
\definecolor{rltgreen}{rgb}{0,0.5,0}
\definecolor{oneblue}{rgb}{0,0,0.75}
\definecolor{brown}{rgb}{0.64,0.16,0.16}
\definecolor{forestgreen}{rgb}{0.13,0.54,0.13}
\definecolor{purple}{rgb}{0.62,0.12,0.94}
\definecolor{dockerblue}{rgb}{0.11,0.56,0.98}
\definecolor{freeblue}{rgb}{0.25,0.41,0.88}
\definecolor{myblue}{rgb}{0,0.2,0.4}
\definecolor{gris}{gray}{0.65}
\begin{document}

\title{Poisson loglinear modeling with linear constraints on the expected cell
frequencies\thanks{This work was partially supported by Grant MTM2009-10072.}}
\author{Nirian Mart\'{\i}n\\Dep. Statistics, Carlos III University of Madrid
\and Leandro Pardo\\Dep. Statistics and O.R., Complutense University of Madrid}
\date{August 6, 2010}
\maketitle

\begin{abstract}
In this paper we consider Poisson loglinear models with linear constraints
(LMLC) on the expected table counts. Multinomial and product multinomial
loglinear models can be obtained by considering that some marginal totals
(linear constraints on the expected table counts) have been prefixed in a
Poisson loglinear model. Therefore with the theory developed in this paper,
multinomial and product multinomial loglinear models can be considered as a
particular case. To carry out inferences on the parameters in the LMLC an
information-theoretic approach is followed from which the classical maximum
likelihood estimators and Pearson chi-square statistics for goodness-of fit
are obtained. In addition, nested hypotheses are proposed as a general
procedure for hypothesis testing. Through a simulation study the
appropriateness of proposed inference tools is illustrated.

\end{abstract}%

\pagestyle{myheadings}%
\markboth{Nirian Mart\'{\i}n and Leandro Pardo}%
{Poisson loglinear modeling with linear constraints}%
\pagenumbering{arabic}%


\section{Introduction\label{Sec1}}

We consider a contingency table with $k$ cells $\boldsymbol{n}=(n_{1}%
,...,n_{k})^{T}$, with $n_{i}$ being the observed frequency associated with
the $i$-th cell ($i=1,...,k$), its distribution is given by a Poisson random
variable and since all of them are mutually independent the joint distribution
of the contingency table is totally specified. Through a loglinear model
$\log\boldsymbol{m}(\boldsymbol{\theta}\mathbf{)=}\boldsymbol{X}%
\boldsymbol{\theta}$ a pattern is established for the mean vector of the
contingency table, $\boldsymbol{m}{(}\boldsymbol{\theta})\equiv(m_{1}%
(\boldsymbol{\theta}),...,m_{k}(\boldsymbol{\theta}))^{T}$, $m_{i}%
(\boldsymbol{\theta})=E[n_{i}]$, $i=1,...,k$, where $\boldsymbol{X}$\ is a
known $k\times t$\ full rank design matrix such that $t\leq k$ and
$\boldsymbol{\theta}=(\theta_{1},...,\theta_{t})^{T}\in\mathbb{R}^{t}$ is the
vector of unknown parameters of the loglinear model.

Let%
\[
\mathcal{C}(\boldsymbol{X})\equiv\{\log\boldsymbol{m}(\boldsymbol{\theta
}\mathbf{):}\log\boldsymbol{m}(\boldsymbol{\theta}\mathbf{)=}\boldsymbol{X}%
\boldsymbol{\theta};\boldsymbol{\theta}\in%
\mathbb{R}
^{t}\}
\]
be the range of loglinear models associated with $\boldsymbol{X}$. We can
observe that $\mathcal{C}(\boldsymbol{X})$\ is the column space of matrix
$\boldsymbol{X}$. A usual convention for loglinear models is to assume that
the vector of $1$'s, $\boldsymbol{J}_{k}\equiv{(1,...,1)}^{T}$, belongs to
$\mathcal{C}(\boldsymbol{X})$, and therefore if a first column $\boldsymbol{J}%
_{k}$\ for $\boldsymbol{X}$\ is considered, the first term $\theta_{1}$ of
$\boldsymbol{\theta}$ is referred to the independent term of the model.

In order to make statistical inference in the class\ of loglinear models
$\mathcal{C}(\boldsymbol{X})$, Cressie and Pardo
\cite{CressiPardo0,CressiePardo} considered for the first time in loglinear
models, minimum $\phi$-divergence estimators and $\phi$-divergence
test-statistics. Later Martin and Pardo \cite{MartinPardo} presented a unified
study for the three different sampling plans (multinomial, product-multinomial
and Poisson).

To study some real situations on the basis of loglinear models, it is
necessary to consider, in addition to loglinear models, some linear
constraints. Loglinear models with linear constraints (LMLC) (see Definition
\ref{Def1}) and product-multinomial sampling were considered for the first
time by Haber and Brown \cite{HaberBrown}. One purpose of this paper is to
consider divergence measures in order to make statistical inference
(estimation and testing) in the class of LMLC but not only with
product-multinomial sampling. We shall present a joint study for different
sampling plans (multinomial, product-multinomial and Poisson). In addition,
this article highlights the fact that the choice of additional lineal
constraints is another way for nesting LMLC, in contrast to the traditional
manner of nesting only log-lineal constraints by reducing the number of
columns of the design matrix. From this idea arises a new way for comparing
LMLC that has not been previously considered in any paper and covers the
preexistent hypothesis-testing techniques as special case (this point will be
clarified in Section \ref{Sec4}).

This article is organized as follows. In Section \ref{Sec1b} we shall consider
some notation as well as some preliminary concepts that will be important in
the other sections of the paper. We pay special attention to the definition of
phi-divergence measures between two non-negative vectors. Section \ref{Sec2}
is devoted to define and study the minimum phi-divergence estimator of LMLC.
The performed \textit{constrained estimation method} will allow us\ to retain
the advantage of dealing with Poisson loglinear models (specially to become
estimation theory easier), even we could have, in fact, a multinomial or
product-multinomial sampling plan. Moreover by an extension of such a method,
if a marginal modeling itself is required, a compact estimation methodology is
provided. As generalization of the constrained maximum likelihood estimation
method, the constrained minimum $\phi$-divergence estimation theory for\ LMLC
is provided. Based also on $\phi$-divergences, in Section \ref{Sec34}%
\ some\ test statistics for LMLC are proposed, specifically Section \ref{Sec3}
is devoted to the problem of goodness-of-fit in LMLC and in Section \ref{Sec4}
the problem of nested hypothesis in LMLC is studied. For both problems the
asymptotic distribution of the $\phi$-divergence\ test statistics under the
null hypothesis are obtained. From such $\phi$-divergence test statistics, in
the case of the goodness-of-fit of LMLC, the classical likelihood ratio and
Pearson chi-square test-statistics, presented in Haber and Brown
\cite{HaberBrown} for multinomial and product-multinomial sampling
schemes,\ are obtained as special case. In Section \ref{Sec5} three hypothesis
tests, which share the aim for testing essentially a marginal model, are
presented. The common framework of the LMLC developed in the previous cited
sections, will allow us to carry out an easier comparison between them. An
example of the potential versatility of such models\ will be shown in Remark
\ref{Rem1}, by considering apparently so different models, such as loglinear
models and marginal models, within the same type of models. Some particular
cases of loglinear models (symmetry, quasi-symmetry and ordinal
quasi-symmetry) on one hand, and a marginal model on the other hand (marginal
homogeneity model) are compared and the exact size and power of their
hypothesis testing is analyzed.

\section{Basic notation and definitions \label{Sec1b}}

By single index notation of $\boldsymbol{n}$ we are able to unify a broad
class of contingency tables, and by convention the terms of multiway
contingency tables can be considered to be located in lexicographical order
but by assigning a single index. For example in the usual double
index\ notation for a a two-way $I\times J$ contingency table $\boldsymbol{n}%
=(n_{11},\allowbreak n_{12},\allowbreak...,\allowbreak n_{1J},...,n_{I1}%
,\allowbreak n_{I2},\allowbreak...,n_{IJ})^{T}$, $n_{ab}$ can be expressed by
a one-to-one index transformation $i=(a-1)J+b$, and therefore $k=IJ$. In a
three-way $I\times J\times K$ contingency table, $\boldsymbol{n}%
=(n_{111},\allowbreak n_{112},\allowbreak...,n_{11K},\allowbreak
...,n_{1J1},\allowbreak n_{1J2},\allowbreak...,n_{1JK},\allowbreak
...,n_{I11},\allowbreak n_{I12},\allowbreak...,n_{I1K},\allowbreak
...,n_{IJ1},\allowbreak n_{IJ2},\allowbreak...,n_{IJK})^{T}$, $n_{abc}$ can be
expressed by a one-to-one index transformation $i=(a-1)JK+(b-1)K+c$, and
therefore $k=IJK$. These models with single index notation are the so-called
coordinate-free models (see Zelterman \cite[Chapter 5]{Zelterman}).

Product-multinomial sampling plan can be considered from a Poisson sampling
plan with some additional linear constraints on the expected cell frequencies.
Since $c$ independent contingency subtables $\boldsymbol{n}_{h}=(n_{h1}%
,...,n_{hk_{h}})^{T}$ are considered in a product-multinomial sampling plan,
the whole contingency table is $\boldsymbol{n}\equiv(\boldsymbol{n}_{1}%
^{T},...,\boldsymbol{n}_{c}^{T})^{T}$\ and $k$ is the summation of the number
of cells in each subtable $k_{h}$, that is $k\equiv\sum_{h=1}^{c}k_{h}$. The
marginal total in each subtable $\sum_{i=1}^{k_{h}}n_{hi}=\boldsymbol{J}%
_{k_{h}}^{T}\boldsymbol{n}_{h}$, $h=1,...,c$\ is prefixed to be $N_{h}\in%
\mathbb{N}
$, and hence the mean vector $\boldsymbol{m}{(}\boldsymbol{\theta
})=(\boldsymbol{m}_{1}(\boldsymbol{\theta}\mathbf{)},...,\boldsymbol{m}%
_{c}(\boldsymbol{\theta}\mathbf{)})^{T}$ with an underlying Poisson sampling
plan $c$ linear constraints are verified%
\begin{equation}
\boldsymbol{J}_{k_{h}}^{T}\boldsymbol{m}_{h}(\boldsymbol{\theta}%
\mathbf{)=}\boldsymbol{J}_{k_{h}}^{T}\boldsymbol{n}_{h}%
,h=1,...,c,\text{\ or\ }\left(
{\displaystyle\bigoplus\limits_{h=1}^{c}}
\boldsymbol{J}_{k_{h}}^{T}\right)  \!\boldsymbol{m}(\boldsymbol{\theta
}\mathbf{)=}\left(
{\displaystyle\bigoplus\limits_{h=1}^{c}}
\boldsymbol{J}_{k_{h}}^{T}\right)  \!\boldsymbol{n},\label{eq1}%
\end{equation}
with $%
{\textstyle\bigoplus\nolimits_{h=1}^{d}}
\boldsymbol{A}_{d}\equiv\operatorname*{diag}\{\boldsymbol{A}_{_{1}%
},...,\boldsymbol{A}_{d}\}$ representing the direct sum of $d$ matrices. In
particular, for multinomial sampling by taking $c=1$ we have%
\begin{equation}
\boldsymbol{J}_{k}^{T}\boldsymbol{m}(\boldsymbol{\theta}\mathbf{)=}%
\boldsymbol{J}_{k}^{T}\boldsymbol{n}.\label{eq2}%
\end{equation}
In what follows $c=0$, i.e. the case where there is no any linear restriction
associated with the sampling plan, will represent that the Poisson sampling
itself is being taken into account.

It is well-known that there are some equivalences between the inferential
results associated with the parameters for the three sampling plans (see for
instance or instance in Lang \cite{Lang,Lang2} and Agresti \cite[Section
14.4]{Agresti}). The main reason why Poisson loglinear model is simpler to
handle is based on the independence of the components of the sampling data.

The parameter space is given by%
\begin{equation}
\Theta=\{\boldsymbol{\theta}\in%
\mathbb{R}
^{t}:\boldsymbol{X}_{0}^{T}\boldsymbol{m}(\boldsymbol{\theta}\mathbf{)=}%
\boldsymbol{X}_{0}^{T}\boldsymbol{n}\},\label{eq3}%
\end{equation}
where $\boldsymbol{X}_{0}\equiv%
{\textstyle\bigoplus\nolimits_{h=1}^{c}}
\boldsymbol{J}_{k_{h}}$ if $c\geq1$ and $\boldsymbol{X}_{0}$ is a vector of
zeros $\boldsymbol{0}_{k}$ if $c=0$ (i.e., $\Theta=%
\mathbb{R}
^{t}$). When $c\geq2$\ a stronger assumption than $\boldsymbol{J}_{k}%
\in\mathcal{C}(\boldsymbol{X})$ is taken into account, $\mathcal{C}%
(\boldsymbol{X}_{0})\subset\mathcal{C}(\boldsymbol{X})$, and therefore if the
first $c$ columns $\boldsymbol{X}_{0}$\ for $\boldsymbol{X}$\ are considered
the $h$-th term $\theta_{h}$ ($h=1,..,c$) of $\boldsymbol{\theta}$ is referred
to the independent term for the model focussed only on the $h$-th contingency subtable.

Haber and Brown \cite{HaberBrown} considered multinomial and product
multinomial LMLC but they did not consider the problem with Poisson sampling.
Definition \ref{Def1} is an extension of the definition given by Haber and
Brown in which Poisson Loglinear models are included.

\begin{Definition}
\label{Def1}In addition of $c$ linear constraints of the sampling scheme,
consider $r\leq t-c$ linear constraints, $\boldsymbol{C}^{T}\boldsymbol{m}%
(\boldsymbol{\theta}\mathbf{)}=\boldsymbol{d}^{\ast}$ , i.e. $\boldsymbol{C} $
and $\boldsymbol{d}^{\ast}$\ are $k\times r$\ and $r\times1$\ matrices
respectively. Once a loglinear model is established through a design matrix
$\boldsymbol{X}$, a loglinear model with linear constraints is a simultaneous
modeling of $\boldsymbol{m}(\boldsymbol{\theta}\mathbf{)}$\ through a
loglinear pattern on one hand and a linear pattern on the other hand%
\begin{equation}
\log\boldsymbol{m}(\boldsymbol{\theta}\mathbf{)=}\boldsymbol{X}%
\boldsymbol{\theta}\quad\text{and}\quad\boldsymbol{L}^{T}\boldsymbol{m}%
(\boldsymbol{\theta}\mathbf{)}=\boldsymbol{d},\label{eq7}%
\end{equation}
being $\boldsymbol{L}{=(}\boldsymbol{X}{_{0},}\boldsymbol{C}{)}$,
$\boldsymbol{d}=(\boldsymbol{n}^{T}\boldsymbol{X}_{0},(\boldsymbol{d}^{\ast
})^{T})^{T}$, for $c\geq1$, and $\boldsymbol{L}{=}\boldsymbol{C}$,
$\boldsymbol{d}=\boldsymbol{d}^{\ast}$, for $c=0$. It is also assumed to hold
$k\geq t-c-r$, and $\operatorname*{rank}(\boldsymbol{L})=\operatorname*{rank}%
(\boldsymbol{L}{,}\boldsymbol{d})=c+r$.
\end{Definition}

Several examples are shown in Haber and Brown \cite{HaberBrown} for $c\geq1$
and an application for $c=0$ and $r\geq1$\ is suggested in Gail \cite[Section
5]{Gail} (for more details see Pardo and Mart\'{\i}n \cite{PardoMartin}).

The parameter space of (\ref{eq7}) is given by%
\begin{equation}
\Theta=\{\boldsymbol{\theta}\mathbf{\in%
\mathbb{R}
}^{t}:\boldsymbol{L}^{T}\boldsymbol{m}(\boldsymbol{\theta}\mathbf{)}%
=\boldsymbol{d}{\}}.\label{eq8}%
\end{equation}
It has been pointed-out that in most practical cases $\boldsymbol{d}^{\ast
}=\boldsymbol{0}_{r}$, actually it holds in all examples of Haber and Brown
\cite{HaberBrown}. Observe that $\boldsymbol{d}^{\ast}=({d}_{1}^{\ast}%
,...,{d}_{r}^{\ast})^{T}$\ has been assumed to be constant, in fact if
${d}_{j}^{\ast}$\ ($j\in\{1,...,r\}$) is proportional to $N\equiv%
{\textstyle\sum\nolimits_{i=1}^{k}}
m_{i}(\boldsymbol{\theta})$ there exists another equivalent constraint where
${d}_{j}^{\ast}={0}$.

In establishing asymptotic properties, we let $N$ tend to infinity, and in
this condition it is assumed that the normalized vector $\boldsymbol{m}^{\ast
}(\boldsymbol{\theta}\mathbf{)=}\boldsymbol{m}(\boldsymbol{\theta}%
\mathbf{)}/N$ remains fixed. For $c\geq1$\ this implies that, as
$N\rightarrow\infty$, the probabilities in each cell remain fixed and
$N_{h}/N$, $h=1,...,c$, remain also fixed.

\begin{Remark}
\label{Rem1}It is interesting to observe that we can consider two cases of
LMLC: \newline\textbf{i)} The classical \textsl{loglinear models }without
linear constraints, which are only defined through the loglinear pattern
$\log\boldsymbol{m}(\boldsymbol{\theta}\mathbf{)=}\boldsymbol{X}%
\boldsymbol{\theta}$ and thus $r=0$ and $\boldsymbol{L}{=}\boldsymbol{X}_{0}%
$.\newline\textbf{ii)} The \textsl{marginal models}, which are only defined
through the linear pattern $\boldsymbol{L}^{T}\boldsymbol{m}%
(\boldsymbol{\theta}\mathbf{)}=\boldsymbol{d}$\ and thus by considering that
$\boldsymbol{X}$ is given by the identity matrix of order $k$, $\boldsymbol{I}%
_{k}$\ (i.e., $k=t$) the loglinear pattern is not itself a restriction.
\end{Remark}

In the particular case $c=0$ and $r=0$, i.e., Poisson loglinear models without
linear constraints, Cressie and Pardo \cite{CressiePardo}\ considered the
problem of testing using divergence measures between probability vectors and
solving the problem of estimation using the maximum likelihood estimator.
Later in Martin and Pardo \cite{MartinPardo} the problem of estimation and
testing was considered using divergence measures between nonnegative vectors
but only for $r=0$. Now in this paper the results obtained in Martin and Pardo
\cite{MartinPardo} are extended for any $r\geq0.$ In this extension we
consider the $\phi$-divergence measure between nonnegative vectors.

Let $\Phi$ be the class of all convex and differentiable functions
$\phi:\left[  0,\infty\right)  \rightarrow%
\mathbb{R}
\cup\left\{  \infty\right\}  $, such that at $x=1$, $\phi\left(  1\right)
=\phi^{\prime}\left(  1\right)  =0$, $\phi^{\prime\prime}\left(  1\right)
>0$. A $\phi$-divergence measure between the $%
\mathbb{R}
_{+}^{k}$-vectors $\boldsymbol{a}=(a_{1},...,a_{k})^{T}$ and $\boldsymbol{b}%
=(b_{1},...,b_{k})^{T}$ is given by%
\begin{equation}
D_{\phi}\left(  \boldsymbol{a}{,}\boldsymbol{b}\right)  =%
{\displaystyle\sum\limits_{i=1}^{k}}
b_{i}\phi\left(  \frac{a_{i}}{b_{i}}\right)  ,\quad\phi\in\Phi,\label{eq4}%
\end{equation}
where $0\phi\left(  0/0\right)  \equiv0$ and $0\phi\left(  p/0\right)  \equiv
p\lim\nolimits_{u\rightarrow\infty}\phi\left(  u\right)  /u$\ conventions are
assumed. These measures cover the traditional ones for probabilistic
arguments, analyzed in Pardo \cite{Pardo}, and all of them share similar
properties. In particular by taking $\lambda\in%
\mathbb{R}
$\ and%
\begin{equation}
\phi_{(\lambda)}(x)=\frac{x^{\lambda+1}-x-\lambda\left(  x-1\right)  }%
{\lambda\left(  \lambda+1\right)  }\text{,}\quad\text{if }\lambda
(\lambda+1)\neq0\text{,}\label{eq5}%
\end{equation}
and $\phi_{(\lambda^{\ast})}(x)=\lim_{\lambda\rightarrow\lambda^{\ast}}%
\phi_{(\lambda)}(x)$, if $\lambda^{\ast}\in\{0,-1\}$, power divergence
measures, introduced in Cressie and Read \cite{CressieRead}, are obtained. The
so-called Kullback divergence measure is obtained through $\phi_{(0)}(x)=x\log
x-x+1$,%
\begin{equation}
D_{Kull}\left(  \boldsymbol{a}{,}\boldsymbol{b}\right)  \equiv D_{\phi_{(0)}%
}\left(  \boldsymbol{a}{,}\boldsymbol{b}\right)  =%
{\displaystyle\sum\limits_{i=1}^{k}}
a_{i}\log\left(  \frac{a_{i}}{b_{i}}\right)  -%
{\displaystyle\sum\limits_{i=1}^{k}}
a_{i}+%
{\displaystyle\sum\limits_{i=1}^{k}}
b_{i},\label{eq6}%
\end{equation}
which was given between two non-negative vectors for the first time in
Brockett \cite{Brockett}. It should be pointed out that the way in which
asymptotic results were obtained in \cite{MartinPardo} is primarily focussed
on the parameter vector, being the mean vector a secondary aim, and therefore
this way is just opposite to the one followed for other works related to
loglinear modeling (see for instance Lang \cite{Lang0}) where the primary aim
is the mean vector itself. These measures and also the methodology for
developing asymptotic results will remain being useful for obtaining the
asymptotic results associated with LMLC. Taking into account Remark
\ref{Rem1}, it is important to clarify that apart from the possibility of
reproducing all inferential results obtained previously in Mart\'{\i}n and
Pardo \cite{MartinPardo}, the new results of this paper are important because
the LMLC cover a broad range of models.

\section{Minimum $\phi$-divergence estimator\label{Sec2}}

The maximum likelihood estimator (MLE) $\widehat{\boldsymbol{\theta}}$ of the
parameter in (\ref{eq7}) can be obtained by maximizing the kernel of the
Poisson log-likelihood%
\[
\ell(\boldsymbol{n}\mathbf{,}\boldsymbol{m}(\boldsymbol{\theta}\mathbf{))}%
\equiv%
{\displaystyle\sum\limits_{i=1}^{k}}
n_{i}\log m_{i}(\boldsymbol{\theta}\mathbf{)}-%
{\displaystyle\sum\limits_{i=1}^{k}}
m_{i}(\boldsymbol{\theta}\mathbf{)},
\]
subject to the constraints $\boldsymbol{L}^{T}\boldsymbol{m}%
(\boldsymbol{\theta}\mathbf{)}=\boldsymbol{d}$, i.e. on the basis of
(\ref{eq8})%
\[
\widehat{\boldsymbol{\theta}}=\arg\max_{\boldsymbol{\theta}\mathbf{\in}\Theta
}\ell(\boldsymbol{n}\mathbf{,}\boldsymbol{m}(\boldsymbol{\theta}\mathbf{))}.
\]
Observe that according to Definition \ref{Def1}, if $\boldsymbol{X}_{0}=%
{\textstyle\bigoplus\nolimits_{h=1}^{c}}
\boldsymbol{J}_{k_{h}}$, which takes part in $\boldsymbol{L}$ as submatrix,
the underlying sampling plan is product-multinomial (or multinomial, if
$c=1$). In what follows even sometimes (product) multinomial sampling will not
be explicitly mentioned, in all results this sampling plan will be covered.

On the basis of (\ref{eq6}) we have%
\[
D_{Kull}(\boldsymbol{n}{,}\boldsymbol{m}(\boldsymbol{\theta}\mathbf{))}=%
{\displaystyle\sum\limits_{i=1}^{k}}
n_{i}\log n_{i}\mathbf{-}%
{\displaystyle\sum\limits_{i=1}^{k}}
n_{i}-\ell(\boldsymbol{n}\mathbf{,}\boldsymbol{m}(\boldsymbol{\theta
}\mathbf{))},
\]
and it is possible also define the MLE of the parameter in (\ref{eq7}) by%
\[
\widehat{\boldsymbol{\theta}}=\arg\min_{\boldsymbol{\theta}\mathbf{\in}\Theta
}D_{Kull}(\boldsymbol{n}{,}\boldsymbol{m}(\boldsymbol{\theta}\mathbf{)).}%
\]
Rather than using MLE, one could use divergence based methods for estimating
the parameters of the loglinear models with linear constraints. On the basis
of (\ref{eq4}) a minimum $\phi$-divergence estimator (M$\phi$E) for a LMLC,
given in Definition \ref{Def1}, is defined as follows.

\begin{Definition}
\label{Def2}For a LMLC (\ref{eq7}) with parameter space (\ref{eq8}), the
M$\phi$E\ is given by%
\begin{equation}
\widehat{\boldsymbol{\theta}}^{\phi}=\arg\min_{\boldsymbol{\theta}\mathbf{\in
}\Theta}D_{\phi}(\boldsymbol{n}{,}\boldsymbol{m}(\boldsymbol{\theta
}\mathbf{))},\label{eq9}%
\end{equation}
with $D_{\phi}(\boldsymbol{n}{,}\boldsymbol{m}(\boldsymbol{\theta}\mathbf{))}%
$\ defined by (\ref{eq4}).
\end{Definition}

In Aitchison and Silvey \cite{AitchisonSylvey} a method for finding MLE's
subject to constraints and its asymptotic distribution theory was developed
for the first time using the Lagrange multiplier method. In Pardo et al.
\cite{PardoPardoZografos} a M$\phi$E procedure for multinomial models was
introduced in which the probabilities depend on unknown parameters that
satisfy some functional relationships. Following the last method but more
generally in the sense that the probabilities are replaced by means, in the
following theorem we present the key result for developing the asymptotic
distribution theory for LMLC, the decomposition of the M$\phi$E for the
parameter vector.

\begin{theorem}
\label{Th1}Suppose that the data $\boldsymbol{n}=\left(  n_{1},...,n_{k}%
\right)  ^{T}$\ are Poisson distributed whose mean vector is given by a LMLC
(\ref{eq7}). Choosing a function $\phi\in\Phi,$ where $\Phi$ was defined in
Section \ref{Sec1}, we have%
\[
\widehat{\boldsymbol{\theta}}^{\phi}=\boldsymbol{\theta}_{0}+\boldsymbol{H}%
(\boldsymbol{\theta}_{0})\boldsymbol{X}^{T}\left(  \frac{\boldsymbol{n}}%
{N}-\boldsymbol{m}^{\ast}{(}\boldsymbol{\theta}_{0})\right)  +o\left(
\left\Vert \frac{\boldsymbol{n}}{N}-\boldsymbol{m}^{\ast}{(}\boldsymbol{\theta
}_{0})\right\Vert \right)  ,
\]
where%
\begin{align*}
\boldsymbol{H}(\boldsymbol{\theta}_{0}) &  \equiv\boldsymbol{I}_{\!\mathcal{F}%
}{(}\boldsymbol{\theta}{_{0})}^{-1}-\boldsymbol{I}_{\!\mathcal{F}}%
{(}\boldsymbol{\theta}{_{0})}^{-1}\boldsymbol{B}{(}\boldsymbol{\theta}%
{_{0}){\left(  \boldsymbol{B}{(}\boldsymbol{\theta}{_{0})}^{T}\boldsymbol{I}%
_{\!\mathcal{F}}{(}\boldsymbol{\theta}{_{0})}^{-1}\boldsymbol{B}%
{(}\boldsymbol{\theta}{_{0})}\right)  }^{-1}}\\
&  \times\boldsymbol{B}{(}\boldsymbol{\theta}{{_{0})}^{T}}\boldsymbol{I}%
{_{\!\mathcal{F}}{(}\boldsymbol{\theta}{_{0})}^{-1},}\\
\boldsymbol{I}_{\!\mathcal{F}}{(}\boldsymbol{\theta}{_{0})} &  {\equiv
}\boldsymbol{X}^{T}\boldsymbol{D}_{\boldsymbol{m}^{\ast}(\boldsymbol{\theta
}_{0})}\boldsymbol{X}\text{ (Fisher information matrix associated with the}\\
&  \text{Poisson loglinear model),}\\
\boldsymbol{B}{(}\boldsymbol{\theta}{_{0})} &  {\equiv}\boldsymbol{X}%
^{T}\boldsymbol{D}_{\boldsymbol{m}^{\ast}(\boldsymbol{\theta}_{0}\mathbf{)}%
}\boldsymbol{L}{,}%
\end{align*}
$\boldsymbol{D}_{\boldsymbol{m}^{\ast}(\boldsymbol{\theta}_{0})}$ is the
diagonal matrix of the normalized vector $\boldsymbol{m}^{\ast}%
(\boldsymbol{\theta}\mathbf{)}${\ and }$\boldsymbol{\theta}_{0}\mathbf{\in
}\Theta$\ is the true and unknown value of the parameter{.}
\end{theorem}

\begin{proof}
We omit the proof because its steps are similar to ones given in Mart\'{\i}n
and Pardo \cite{MartinPardo0} with the differences motivated because in the
cited paper only multinomial sampling was considered. $\square$
\end{proof}

In the next theorem we obtain the asymptotic distribution of $\widehat
{\boldsymbol{\theta}}^{\phi}$\ as well as of $\boldsymbol{m}(\widehat
{\boldsymbol{\theta}}^{\phi})$.

\begin{theorem}
\label{Th2}Suppose that the data $\boldsymbol{n}=\left(  n_{1},...,n_{k}%
\right)  ^{T}$\ are Poisson distributed whose mean vector is given by a LMLC
(\ref{eq7}). Choosing a function $\phi\in\Phi,$ where $\Phi$ was defined in
Section \ref{Sec1}, we have\newline a)%
\begin{equation}
\sqrt{N}(\widehat{\boldsymbol{\theta}}^{\phi}-\boldsymbol{\theta}_{0}%
)\overset{\mathcal{L}}{\underset{N\rightarrow\infty}{\longrightarrow}%
}\mathcal{N}(\boldsymbol{0}_{{t}},\boldsymbol{H}{(}\boldsymbol{\theta}%
{_{0}\mathbf{))}}\label{eq10}%
\end{equation}
where $\boldsymbol{H}{(}\boldsymbol{\theta}{_{0}\mathbf{)}}$ is defined in
Theorem \ref{Th1}, \textquotedblleft$\overset{\mathcal{L}}{\underset
{N\rightarrow\infty}{\longrightarrow}}$\textquotedblright\ denotes convergence
in law (or distribution) and\newline b)%
\begin{equation}
\frac{1}{\sqrt{N}}(\boldsymbol{m}(\widehat{\boldsymbol{\theta}}^{\phi
})-\boldsymbol{m}(\boldsymbol{\theta}_{0}))\overset{\mathcal{L}}%
{\underset{N\rightarrow\infty}{\longrightarrow}}\mathcal{N}(\boldsymbol{0}%
_{k},\boldsymbol{\Sigma})\label{eq11}%
\end{equation}
where%
\[
\boldsymbol{\Sigma}{\equiv}\boldsymbol{D}_{\boldsymbol{m}^{\ast}%
(\boldsymbol{\theta}_{0}\mathbf{)}}\boldsymbol{XH}{(}\boldsymbol{\theta}%
{_{0})}\boldsymbol{X}^{T}\boldsymbol{D}_{\boldsymbol{m}^{\ast}%
(\boldsymbol{\theta}_{0}\mathbf{)}}.
\]

\end{theorem}

\begin{proof}
Result a) follows by Theorem \ref{Th1}\ and taking into account (see Haberman
\cite{Haberman})%
\begin{equation}
\frac{1}{\sqrt{N}}(\boldsymbol{n}-\boldsymbol{m}(\boldsymbol{\theta}%
_{0}))\overset{\mathcal{L}}{\underset{N\rightarrow\infty}{\longrightarrow}%
}\mathcal{N}(\boldsymbol{0}_{k},\boldsymbol{D}_{\boldsymbol{m}^{\ast
}(\boldsymbol{\theta}_{0})}).\label{eq12}%
\end{equation}
Part b) follows by a) and applying delta method (see for instance Agresti
\cite[Sections 14.1.2, 14.1.3]{Agresti}). $\square$
\end{proof}

In the next theorem a result related to a simplification of the expression of
the asymptotic variance-covariance matrices of Theorem \ref{Th2} is shown.

\begin{theorem}
\label{Th3}When%
\begin{equation}
\boldsymbol{X}{=(}\boldsymbol{L}{,}\boldsymbol{W}{)}\label{eq13}%
\end{equation}
we have%
\[
\boldsymbol{H}{(}\boldsymbol{\theta}{_{0})=}(\boldsymbol{X}^{T}\boldsymbol{D}%
_{\boldsymbol{m}^{\ast}(\boldsymbol{\theta}_{0})}\boldsymbol{X}{)}%
^{-1}-(\boldsymbol{L}^{T}\boldsymbol{D}_{\boldsymbol{m}^{\ast}%
(\boldsymbol{\theta}_{0})}\boldsymbol{L}{)}^{-1}\oplus\boldsymbol{0}%
_{(t-c-r)\times(t-c-r)}%
\]
and%
\[
\boldsymbol{\Sigma}{=}\boldsymbol{D}_{\boldsymbol{m}^{\ast}(\boldsymbol{\theta
}_{0}\mathbf{)}}^{\frac{1}{2}}(\boldsymbol{A}_{X}{(}\boldsymbol{\theta}%
{_{0})-}\boldsymbol{A}_{L}{(}\boldsymbol{\theta}{_{0})})\boldsymbol{D}%
_{\boldsymbol{m}^{\ast}(\boldsymbol{\theta}_{0}\mathbf{)}}^{\frac{1}{2}},
\]
where%
\begin{align*}
\boldsymbol{A}_{X}{(}\boldsymbol{\theta}{_{0})} &  {\equiv}\boldsymbol{D}%
{{_{\boldsymbol{m}^{\ast}(\boldsymbol{\theta}_{0})}^{\frac{1}{2}}}%
}\boldsymbol{X}{\left(  \boldsymbol{X}^{T}\boldsymbol{D}_{\boldsymbol{m}%
^{\ast}(\boldsymbol{\theta}_{0})}\boldsymbol{X}\right)  ^{-1}}\boldsymbol{X}%
{^{T}}\boldsymbol{D}{_{\boldsymbol{m}^{\ast}(\boldsymbol{\theta}_{0})}%
^{\frac{1}{2}},}\\
\boldsymbol{A}_{L}{(}\boldsymbol{\theta}{_{0})} &  {\equiv}\boldsymbol{D}%
{{_{\boldsymbol{m}^{\ast}(\boldsymbol{\theta}_{0})}^{\frac{1}{2}}}%
}\boldsymbol{L}{\left(  \boldsymbol{L}^{T}\boldsymbol{D}_{\boldsymbol{m}%
^{\ast}(\boldsymbol{\theta}_{0})}\boldsymbol{L}\right)  ^{-1}}\boldsymbol{L}%
{^{T}}\boldsymbol{D}{_{\boldsymbol{m}^{\ast}(\boldsymbol{\theta}_{0})}%
^{\frac{1}{2}}.}%
\end{align*}

\end{theorem}

We can observe that $\boldsymbol{A}_{X}{(}\boldsymbol{\theta}{_{0})}$ and
$\boldsymbol{A}_{L}{(}\boldsymbol{\theta}{_{0})}$\ are projector matrices
respectively on column spaces $\mathcal{C}(\boldsymbol{D}{{_{{\boldsymbol{m}%
}^{\ast}({\boldsymbol{\theta}}_{0})}^{\frac{1}{2}}}}\boldsymbol{X})$\ and
$\mathcal{C}(\boldsymbol{D}{{_{{\boldsymbol{m}}^{\ast}({\boldsymbol{\theta}%
}_{0})}^{\frac{1}{2}}}}\boldsymbol{L})$.\medskip

\begin{proof}
By starting through a identity matrix,%
\begin{align*}
\boldsymbol{I}_{{t}} &  =\left(  \boldsymbol{X}^{T}\boldsymbol{D}%
_{\boldsymbol{m}^{\ast}(\boldsymbol{\theta}_{0}\mathbf{)}}\boldsymbol{X}%
\right)  ^{-1}\boldsymbol{X}^{T}\boldsymbol{D}_{\boldsymbol{m}^{\ast
}(\boldsymbol{\theta}_{0}\mathbf{)}}\boldsymbol{X}\\
&  =\left(  \boldsymbol{X}^{T}\boldsymbol{D}_{\boldsymbol{m}^{\ast
}(\boldsymbol{\theta}_{0}\mathbf{)}}\boldsymbol{X}\right)  ^{-1}%
\boldsymbol{X}^{T}\boldsymbol{D}_{\boldsymbol{m}^{\ast}(\boldsymbol{\theta
}_{0}\mathbf{)}}\left(  \boldsymbol{L}{,}\boldsymbol{W}\right) \\
&  =\left(  \left(  \boldsymbol{X}^{T}\!\boldsymbol{D}_{\boldsymbol{m}^{\ast
}(\boldsymbol{\theta}_{0}\mathbf{)}}\!\boldsymbol{X}\right)  ^{-1}%
\!\boldsymbol{X}^{T}\!\boldsymbol{D}_{\boldsymbol{m}^{\ast}(\boldsymbol{\theta
}_{0}\mathbf{)}}\!\boldsymbol{L}{,\left(  \boldsymbol{X}^{T}\!\boldsymbol{D}%
_{\boldsymbol{m}^{\ast}(\boldsymbol{\theta}_{0}\mathbf{)}}\!\boldsymbol{X}%
\right)  ^{-1}}\!\boldsymbol{X}{^{T}}\!\boldsymbol{D}{_{\boldsymbol{m}^{\ast
}(\boldsymbol{\theta}_{0}\mathbf{)}}}\!\boldsymbol{W}\right)  ,
\end{align*}
it is obtained that%
\[
\boldsymbol{G}_{t\times(c+r)}=\left(  \boldsymbol{X}^{T}\boldsymbol{D}%
_{\boldsymbol{m}^{\ast}(\boldsymbol{\theta}_{0}\mathbf{)}}\boldsymbol{X}%
\right)  ^{-1}\boldsymbol{X}^{T}\boldsymbol{D}_{\boldsymbol{m}^{\ast
}(\boldsymbol{\theta}_{0}\mathbf{)}}\boldsymbol{L}{=}%
\begin{pmatrix}
\boldsymbol{I}_{c+r}\\
\boldsymbol{0}_{(t-c-r)\times(c+r)}%
\end{pmatrix}
.
\]
Therefore%
\begin{align*}
\boldsymbol{H}{(}\boldsymbol{\theta}{_{0})} &  =\left(  \boldsymbol{X}%
^{T}\boldsymbol{D}_{\boldsymbol{m}^{\ast}(\boldsymbol{\theta}_{0}\mathbf{)}%
}\boldsymbol{X}\right)  ^{-1}-\left(  \boldsymbol{X}^{T}\boldsymbol{D}%
_{\boldsymbol{m}^{\ast}(\boldsymbol{\theta}_{0}\mathbf{)}}\boldsymbol{X}%
\right)  ^{-1}\boldsymbol{X}^{T}\boldsymbol{D}_{\boldsymbol{m}^{\ast
}(\boldsymbol{\theta}_{0}\mathbf{)}}\boldsymbol{L}\\
&  {\times{\left(  \boldsymbol{L}^{T}\boldsymbol{D}_{\boldsymbol{m}^{\ast
}(\boldsymbol{\theta}_{0}\mathbf{)}}\boldsymbol{X}\left(  \boldsymbol{X}%
^{T}\boldsymbol{D}_{\boldsymbol{m}^{\ast}(\boldsymbol{\theta}_{0}\mathbf{)}%
}\boldsymbol{X}\right)  ^{-1}\boldsymbol{X}^{T}\boldsymbol{D}_{\boldsymbol{m}%
^{\ast}(\boldsymbol{\theta}_{0}\mathbf{)}}\boldsymbol{L}\right)  }^{-1}}\\
&  \times\boldsymbol{L}^{T}\boldsymbol{D}_{\boldsymbol{m}^{\ast}%
(\boldsymbol{\theta}_{0}\mathbf{)}}\boldsymbol{X}{\left(  \boldsymbol{X}%
^{T}\boldsymbol{D}_{\boldsymbol{m}^{\ast}(\boldsymbol{\theta}_{0}\mathbf{)}%
}\boldsymbol{X}\right)  ^{-1}}\\
&  =\left(  \boldsymbol{X}^{T}\!\boldsymbol{D}_{\boldsymbol{m}^{\ast
}(\boldsymbol{\theta}_{0}\mathbf{)}}\!\boldsymbol{X}\right)  ^{-1}%
\!-\!\boldsymbol{G}_{t\times(c+r)}\!\left(  \boldsymbol{L}^{T}\!\boldsymbol{D}%
_{\boldsymbol{m}^{\ast}(\boldsymbol{\theta}_{0}\mathbf{)}}\!{{\left(
\boldsymbol{L}{,}\boldsymbol{W}\right)  }}\boldsymbol{G}_{t\times c}%
^{T}\right)  ^{-1}\!\boldsymbol{G}_{t\times(c+r)}^{T}\\
&  =\left(  \boldsymbol{X}^{T}\boldsymbol{D}_{\boldsymbol{m}^{\ast
}(\boldsymbol{\theta}_{0}\mathbf{)}}\boldsymbol{X}\right)  ^{-1}%
-\boldsymbol{G}_{t\times(c+r)}\left(  \boldsymbol{L}^{T}\boldsymbol{D}%
_{\boldsymbol{m}^{\ast}(\boldsymbol{\theta}_{0}\mathbf{)}}\boldsymbol{L}%
\right)  ^{-1}\boldsymbol{G}_{t\times(c+r)}^{T}\\
&  =\left(  \boldsymbol{X}^{T}\boldsymbol{D}_{\boldsymbol{m}^{\ast
}(\boldsymbol{\theta}_{0}\mathbf{)}}\boldsymbol{X}\right)  ^{-1}-\left(
\boldsymbol{L}^{T}\boldsymbol{D}_{\boldsymbol{m}^{\ast}(\boldsymbol{\theta
}_{0})}\boldsymbol{L}\right)  ^{-1}\oplus\boldsymbol{0}_{(t-c-r)\times
(t-c-r)}.
\end{align*}
On the other hand, the expression of $\boldsymbol{\Sigma}$ is obtained
replacing the new expression of $\boldsymbol{H}{(}\boldsymbol{\theta}{_{0})} $
inside its original definition in Theorem \ref{Th2}. $\square$
\end{proof}

\begin{Remark}
\label{Rem2}When $r=0$ and $c=1$ (classical multinomial loglinear model),
$\boldsymbol{X}{=(}\boldsymbol{J}_{k}{,}\boldsymbol{W}{)}$, $\boldsymbol{J}%
_{k}^{T}\boldsymbol{D}_{\boldsymbol{m}^{\ast}(\boldsymbol{\theta}_{0}%
)}\boldsymbol{J}_{k}=1$ and thus it holds $\boldsymbol{H}{(}\boldsymbol{\theta
}{_{0})=\allowbreak}(\boldsymbol{X}^{T}\boldsymbol{D}_{\boldsymbol{m}^{\ast
}(\boldsymbol{\theta}_{0})}\boldsymbol{X}{)}^{-1}-1\oplus\boldsymbol{0}%
_{(t-1)\times(t-1)}$ as direct application of Theorem \ref{Th3}. If we pay
attention to the structure of this variance-covariance matrix, we can observe
that the $\boldsymbol{\theta}=(\theta_{1},\overline{\boldsymbol{\theta}}%
^{T})^{T}$ is partitioned in such a way that once the part associated with
$\boldsymbol{W}$, $\overline{\boldsymbol{\theta}}=(\theta_{2},...,\theta
_{t})^{T}$,\ is known,\ the first component $\theta_{1}$\ can be obtained
through $\overline{\boldsymbol{\theta}}$ and the linear constraint. This is
the reason why in the traditional multinomial loglinear modeling the dimension
of the parameter space is $t-1$ instead of $t$ and $\theta_{1}%
=N/(\boldsymbol{J}_{k}^{T}\exp\{\boldsymbol{W}\overline{\boldsymbol{\theta}%
}\})$ is the redundant component of the multinomial loglinear model. When
$c\geq0$ and $r\geq1$, it is possible to partitionate and interpret any
parameter vector $\boldsymbol{\theta}$\ in terms of (\ref{eq13}). Due to space
limitation, we omit it in a formal way. In a less formal way we can say that
making transformation on the design or restrictions matrices, it is possible
to obtain LMLC with an structure for the design matrix like in (\ref{eq13}).
The part of the parameters associated with matrix $\boldsymbol{W}$, are
\textquotedblleft free parameters\textquotedblright, while the rest of the
terms are determinated through a function. It is frequent to find textbooks
that consider only free parameters for making statistical inferences.
\end{Remark}

\section{$\phi$-divergence test statistics\label{Sec34}}

\subsection{Goodness-of-fit\label{Sec3}}

Classical measures for assessing the goodness-of-fit of categorical data
models, estimated by MLE, are the likelihood ratio test statistic, sometimes
referred to as the deviance statistic,%
\begin{equation}
G^{2}(\widehat{\boldsymbol{\theta}})=2%
{\displaystyle\sum\limits_{{i}=1}^{k}}
\left(  n_{i}\log\frac{n_{i}}{{m}_{{i}}(\widehat{\boldsymbol{\theta}})}%
-(n_{i}-{m}_{{i}}(\widehat{\boldsymbol{\theta}}))\right)  ,\label{eq18}%
\end{equation}
and Pearson chi-square test statistic%
\begin{equation}
{X}^{2}(\widehat{\boldsymbol{\theta}})=%
{\displaystyle\sum\limits_{{i}=1}^{k}}
\frac{(n_{i}-{m}_{{i}}(\widehat{\boldsymbol{\theta}}))^{2}}{{m}_{{i}}%
(\widehat{\boldsymbol{\theta}})}.\label{eq19}%
\end{equation}
In Haber an Brown \cite{HaberBrown} the asymptotic distribution of a classical
goodness-of-fit test-statistics for LMLC when the sampling scheme is (product)
multinomial ($c\geq1$, $r\geq0$) was established. On the other hand in
Mart\'{\i}n and Pardo \cite{MartinPardo} divergence based goodness-of-fit
test-statistics were analyzed for loglinear models under the three sampling
schemes ($c\geq0$) when none constraint additional to the sampling ones are
considered ($r=0$). In this section we extend the previous result to the
important context in which $r>0$. In this framework the family of $\phi
$-divergence test statistics is given by%
\begin{equation}
T^{\phi_{1}}(\widehat{\boldsymbol{\theta}}^{\phi_{2}})=\frac{2}{\phi
_{1}^{\prime\prime}(1)}D_{\phi_{1}}(\boldsymbol{n}{,}\boldsymbol{m}%
(\widehat{\boldsymbol{\theta}}^{\phi_{2}})).\label{eq20}%
\end{equation}
Observe that while the divergence based estimator is associated with a
specific\ $\phi_{2}$\ function, the divergence based test-statistic is
associated with a function $\phi_{1}$, not necessarily equal to $\phi_{2}$, in
fact while $G^{2}(\widehat{\boldsymbol{\theta}})=T^{\phi_{(0)}}(\widehat
{\boldsymbol{\theta}}^{\phi_{(0)}})$\ where $\phi_{(0)}(x)=x\log x-x+1$, it
holds ${X}^{2}(\widehat{\boldsymbol{\theta}})=T^{\phi_{(1)}}(\widehat
{\boldsymbol{\theta}}^{\phi_{(0)}})$ where $\phi_{(1)}(x)=\frac{1}{2}%
(x-1)^{2}$.

In the following theorem we establish that the asymptotic distribution of the
family of $\phi$-divergence test statistics, $T^{\phi_{1}}(\widehat
{\boldsymbol{\theta}}^{\phi_{2}})$, is a chi-square with $k-t+c$ degrees of
freedom ($\chi_{k-t+c}^{2}$). Therefore, we do not accept the null hypothesis
in which the model is said to be (\ref{eq7}) if $T^{\phi_{1}}(\widehat
{\boldsymbol{\theta}}^{\phi_{2}})>c$, where $c$ is specified so that the size
of the test is $\alpha$, $\Pr(T^{\phi_{1}}(\widehat{\boldsymbol{\theta}}%
^{\phi_{2}})>c\mid H_{Null})=\Pr(\chi_{k-t+c}^{2}>c)=\alpha$,\ i.e.
$c\equiv\chi_{k-t+c}^{2}(1-\alpha)$ is the ($1-\alpha$)-th quantile of a
$\chi_{k-t+c}^{2}$\ distribution.

\begin{theorem}
\label{Th7}Suppose that the data $\boldsymbol{n}=\left(  n_{1},...,n_{k}%
\right)  ^{T}$ are Poisson distributed. Choose the function $\phi_{1},\phi
_{2}\in\Phi$, where $\Phi$ was defined in Section \ref{Sec1}. Then, for
testing%
\begin{align}
\quad H_{Null}\!\!\!  & :\log\boldsymbol{m}(\boldsymbol{\theta})\in
\mathcal{C}(\boldsymbol{X})\text{ and }\boldsymbol{\theta}\in\Theta
=\{\boldsymbol{\theta}^{\prime}\in\mathbb{R}^{t}:\boldsymbol{L}^{T}%
\boldsymbol{m}{(}\boldsymbol{\theta}^{\prime}\mathbf{)=}\boldsymbol{d}%
\}\text{,}\label{eq24}\\
H_{Alt}\!\!\!  & :\log\boldsymbol{m}(\boldsymbol{\theta})\notin\mathcal{C}%
(\boldsymbol{X})\text{ or }\boldsymbol{\theta}\notin\Theta
=\{\boldsymbol{\theta}^{\prime}\in\mathbb{R}^{t}:\boldsymbol{L}^{T}%
\boldsymbol{m}{(}\boldsymbol{\theta}^{\prime}\mathbf{)=}\boldsymbol{d}%
\}\text{,}\nonumber
\end{align}
the asymptotic null distribution of the test statistic $T^{\phi_{1}}%
(\widehat{\boldsymbol{\theta}}^{\phi_{2}})$, given in (\ref{eq20}), is
chi-squared with $k-t+c$ degrees of freedom.
\end{theorem}

\begin{proof}
We consider the function $f\left(  x,y\right)  =x\phi_{1}\left(  y/x\right)
$. A second order Taylor's expansion of $f(\frac{n_{{i}}}{N},m_{{i}}^{\ast
}(\widehat{\boldsymbol{\theta}}^{\phi_{2}}))$ about $(m_{{i}}^{\ast
}(\boldsymbol{\theta}_{0}),m_{{i}}^{\ast}(\boldsymbol{\theta}_{0}))$ gives%
\[
f\left(  \frac{n_{{i}}}{N},m_{{i}}^{\ast}(\widehat{\boldsymbol{\theta}}%
^{\phi_{2}})\right)  =\frac{\phi_{1}^{\prime\prime}(1)}{2}\frac{\left(
\dfrac{n_{{i}}}{N}-m_{{i}}^{\ast}(\widehat{\boldsymbol{\theta}}^{\phi_{2}%
})\right)  ^{2}}{m_{{i}}^{\ast}(\boldsymbol{\theta}_{0})}+o_{P}(N^{-1}%
);\quad{i}=1,...,k.
\]
Taking into account%
\[
T^{\phi_{1}}(\widehat{\boldsymbol{\theta}}^{\phi_{2}})=\frac{2N}{\phi
_{1}^{\prime\prime}(1)}\sum_{{i}=1}^{k}f\left(  \frac{n_{{i}}}{N},m_{{i}%
}^{\ast}(\widehat{\boldsymbol{\theta}}^{\phi_{2}})\right)  =\sum_{{i}=1}%
^{k}Z_{{i}}^{2}+{o}_{P}\left(  1\right)  ,
\]
where%
\[
Z_{{i}}\equiv\frac{n_{{i}}-m_{{i}}(\widehat{\boldsymbol{\theta}}^{\phi_{2}}%
)}{\sqrt{m_{{i}}(\boldsymbol{\theta}_{0})}},\text{ }{i}=1,...,k,
\]
we obtain the following vectorial expression
\[
T^{\phi_{1}}(\widehat{\boldsymbol{\theta}}^{\phi_{2}})=\boldsymbol{Z}%
^{T}\boldsymbol{Z}{+o}_{P}\left(  1\right)  ,
\]
being%
\[
\boldsymbol{Z}{=(Z_{{1}},...,Z_{{k}})}^{T}\equiv\boldsymbol{D}_{\boldsymbol{m}%
{(}\boldsymbol{\theta}_{0})}^{-\frac{1}{2}}(\boldsymbol{n}{-}\boldsymbol{m}%
(\widehat{\boldsymbol{\theta}}^{\phi_{2}})).
\]
The random vector $\boldsymbol{Z}$\ is asymptotically normal distributed with
mean vector zero and asymptotic variance-covariance matrix%
\begin{equation}
\boldsymbol{T}^{\ast}\equiv\boldsymbol{I}_{k}{-}\boldsymbol{A}{_{0}{(}%
}\boldsymbol{\theta}{{_{0})}-}\boldsymbol{D}{_{\boldsymbol{m}^{\ast
}(\boldsymbol{\theta}_{0}\mathbf{)}}^{\frac{1}{2}}}\boldsymbol{XH}%
{(}\boldsymbol{\theta}{{_{0})}}\boldsymbol{X}{^{T}}\boldsymbol{D}%
{_{\boldsymbol{m}^{\ast}(\boldsymbol{\theta}_{0}\mathbf{)}}^{\frac{1}{2}}%
,}\label{eq21}%
\end{equation}
where $\boldsymbol{A}{_{0}{(}}\boldsymbol{\theta}{{_{0})}}$\ is given by%
\begin{equation}
\boldsymbol{A}{_{0}{(}}\boldsymbol{\theta}{{_{0})=}}\boldsymbol{D}%
{{_{\boldsymbol{m}^{\ast}(\boldsymbol{\theta}_{0})}^{\frac{1}{2}}}%
}\boldsymbol{X}\text{{}}_{0}{\left(  \boldsymbol{X}\text{{}}_{0}%
^{T}\boldsymbol{D}_{\boldsymbol{m}^{\ast}(\boldsymbol{\theta}_{0}%
)}\boldsymbol{X}\text{{}}_{0}\right)  }^{-1}\boldsymbol{X}{{}_{0}^{T}%
}\boldsymbol{D}{_{\boldsymbol{m}^{\ast}(\boldsymbol{\theta}_{0})}^{\frac{1}%
{2}}.}\label{eq23}%
\end{equation}
Then, the asymptotic distribution of the $\phi$-divergence test statistic
$T^{\phi_{1}}(\widehat{\boldsymbol{\theta}}^{\phi_{2}})$ will be a chi-square
iff the matrix $\boldsymbol{T}^{\ast}$ is idempotent and symmetric. It is
clear that $\boldsymbol{T}^{\ast}$ is symmetric, and to establish that it is
idempotent we have%
\[
(\boldsymbol{T}^{\ast})^{2}=\boldsymbol{SS}{-}\boldsymbol{SK}{-}%
\boldsymbol{KS}{+}\boldsymbol{KK}=\boldsymbol{S}{-}\boldsymbol{K}%
{-}\boldsymbol{K}{+}\boldsymbol{K}=\boldsymbol{T}^{\ast},
\]
where $\boldsymbol{S}{=}\boldsymbol{I}{_{k}-}\boldsymbol{A}{_{0}{(}%
}\boldsymbol{\theta}{{_{0})}}$\ and $\boldsymbol{K}{=}\boldsymbol{D}%
_{\boldsymbol{m}^{\ast}(\boldsymbol{\theta}_{0}\mathbf{)}}^{\frac{1}{2}%
}\boldsymbol{XH}{(}\boldsymbol{\theta}{_{0})}\boldsymbol{X}^{T}\boldsymbol{D}%
_{\boldsymbol{m}^{\ast}(\boldsymbol{\theta}_{0}\mathbf{)}}^{\frac{1}{2}}$. The
degrees of freedom of the chi-squared distributed statistic $T^{\phi_{1}%
}(\widehat{\boldsymbol{\theta}}^{\phi_{2}})$ coincides with the trace of the
matrix $\boldsymbol{T}^{\ast}$, i.e. $k-t+c$. $\square$
\end{proof}

\begin{Remark}
\label{Rem5}When $\mathcal{C}(\boldsymbol{L})\subset\mathcal{C}(\boldsymbol{X}%
)$, because $\mathcal{C(}\boldsymbol{D}_{\boldsymbol{m}^{\ast}%
(\boldsymbol{\theta}_{0}\mathbf{)}}^{\frac{1}{2}}\boldsymbol{L}{)}%
\subset\mathcal{C(}\boldsymbol{D}_{\boldsymbol{m}^{\ast}(\boldsymbol{\theta
}_{0}\mathbf{)}}^{\frac{1}{2}}\boldsymbol{X}{)}$ it holds $\boldsymbol{A}%
_{X}{(}${$\boldsymbol{\theta}$}${{_{0}})}\boldsymbol{D}_{\boldsymbol{m}^{\ast
}(\boldsymbol{\theta}_{0}\mathbf{)}}^{\frac{1}{2}}\boldsymbol{L}%
=\boldsymbol{D}_{\boldsymbol{m}^{\ast}(\boldsymbol{\theta}_{0}\mathbf{)}%
}^{\frac{1}{2}}\boldsymbol{L}$, which means that (\ref{eq21}) is given by%
\begin{equation}
\boldsymbol{T}^{\ast}=\boldsymbol{I}_{k}{-}\boldsymbol{A}{_{0}{(}%
}\boldsymbol{\theta}{{_{0})-}}\boldsymbol{A}_{X}{(}\boldsymbol{\theta}{_{0}%
)+}\boldsymbol{A}_{L}{(}\boldsymbol{\theta}{_{0}).}\label{eq22}%
\end{equation}
From this expression it is concluded that when there is no any sampling
constraint ($r=0\Rightarrow\boldsymbol{A}_{L}{(}\boldsymbol{\theta}{_{0}%
)=}\boldsymbol{A}{_{0}{(}}\boldsymbol{\theta}{{_{0})}}$), the
variance-covariance matrix (\ref{eq21}) of the random vector $\boldsymbol{Z}$,
under the assumption that the model of the null hypothesis in (\ref{eq24})
holds, have a common expression, $\boldsymbol{T}^{\ast}=\boldsymbol{I}_{k}%
{-}\boldsymbol{A}_{X}{(}\boldsymbol{\theta}{_{0})}$, for the three sampling
schemes ($c\geq0$).
\end{Remark}

\subsection{Nested hypothesis \label{Sec4}}

Two models are said to be nested if one of them can be obtained from the other
one as special case. This general definition for linear models (see Chatterjee
and Hadi \cite[page 65]{Chatterjee}) can be applied to two loglinear models,
whose design matrices are given by $\boldsymbol{X}_{1}$\ and $\boldsymbol{X}%
_{2}$, in such a way that the first one is said to be nested within the second
one if $\mathcal{C}(\boldsymbol{X}_{1})\subset\mathcal{C}(\boldsymbol{X}_{2}%
)$. Observe that $\operatorname{rank}(\boldsymbol{X}_{1})=t_{1}\leq
\operatorname{rank}(\boldsymbol{X}_{2})=t_{2}$, and therefore $\Theta
_{1}=\{\boldsymbol{\theta}^{\prime}\in\mathbb{R}^{t_{1}}:\boldsymbol{X}%
_{0}^{T}\boldsymbol{m}{(}\boldsymbol{\theta}^{\prime}\mathbf{)=}%
\boldsymbol{X}_{0}^{T}\boldsymbol{n}\}\subset\Theta_{2}=\{\boldsymbol{\theta
}^{\prime}\in\mathbb{R}^{t_{2}}:\boldsymbol{X}_{0}^{T}\boldsymbol{m}%
{(}\boldsymbol{\theta}^{\prime}\mathbf{)=}\boldsymbol{X}_{0}^{T}%
\boldsymbol{n}\}$, i.e. $\forall${$\boldsymbol{\theta}$}$_{1}\in\Theta_{1}$
$\exists${$\boldsymbol{\theta}$}$_{2}\in\Theta_{2}:\boldsymbol{m}%
(\boldsymbol{\theta}_{1})=\boldsymbol{m}(\boldsymbol{\theta}_{2})$. Moreover,
if $\boldsymbol{X}_{2}=(\boldsymbol{X}_{1},\boldsymbol{Y}_{2})$, where
$\operatorname{rank}(\boldsymbol{Y}_{2})=s_{2}$ (i.e., $t_{2}=t_{1}+s_{2}$),
by considering {$\boldsymbol{\theta}$}$_{2}=(${$\boldsymbol{\theta}$}$_{1}%
^{T},\boldsymbol{0}_{s_{2}}^{T})^{T}$, it holds $\boldsymbol{m}%
(\boldsymbol{\theta}_{1})=\boldsymbol{m}(\boldsymbol{\theta}_{2})$, and thus
the loglinear model defined by $\boldsymbol{X}_{1}$\ is nested within the
loglinear model defined by $\boldsymbol{X}_{2}$. In order to clarify that this
is a particular case of nested model, a loglinear model defined by
$\boldsymbol{X}_{1}$\ is said to be a reduced loglinear model of
$\boldsymbol{X}_{2}=(\boldsymbol{X}_{1},\boldsymbol{Y}_{2})$.

In the following definition we consider a sequence of design matrices
$\{\boldsymbol{X}_{b}\}_{b=1}^{B}$\ so that the loglinear model associated
with $\boldsymbol{X}_{b+1}{=(}\boldsymbol{L}{,}\boldsymbol{W}_{b+1}{)}$\ is a
reduced loglinear model of $\boldsymbol{X}_{b}{=(}\boldsymbol{L}%
{,}\boldsymbol{W}_{b}{)}$, $b=1,...,B-1$, which means that $\boldsymbol{W}%
_{b+1} $ is a submatrix of $\boldsymbol{W}_{b}$. Such matrices define a
sequence of LMLC that share the same linear constraints.

\begin{Definition}
\label{Def6}The sequence of LMLC%
\begin{equation}
\log\boldsymbol{m}(\boldsymbol{\theta}_{b}\mathbf{)=}\boldsymbol{X}%
_{b}\boldsymbol{\theta}_{b}\quad\text{and}\quad\boldsymbol{L}^{T}%
\boldsymbol{m}(\boldsymbol{\theta}_{b}\mathbf{)}=\boldsymbol{d}{,}\label{eq27}%
\end{equation}
where $\boldsymbol{X}_{b}=(\boldsymbol{x}{_{1},...,}\boldsymbol{x}{_{t-b+1}})
$, $b\in\{1,...,t-c-r\}$, is called the $b$-th \textsf{reduced LMLC through
the parameter}, because by reducing one unit the dimension of the parameter
space $\Theta_{b}\equiv\{\boldsymbol{\theta}_{b}\in\mathbf{{%
\mathbb{R}
}}^{t-b+1}:\boldsymbol{L}^{T}\boldsymbol{m}(\boldsymbol{\theta}_{b}%
\mathbf{)=}\boldsymbol{d}\}$, it holds $M_{b+1}\subset M_{b}$\ where%
\[
M_{b}\equiv\{\boldsymbol{m}{(}\boldsymbol{\theta}_{b}\mathbf{)}\in\mathbf{{%
\mathbb{R}
}}^{k}:\log\boldsymbol{m}{(}\boldsymbol{\theta}_{b}\mathbf{)}=\boldsymbol{X}%
_{b}\boldsymbol{\theta}_{b},\boldsymbol{\theta}_{b}\in\Theta_{b}\}.
\]

\end{Definition}

In the following definition we consider a sequence of constraints
$\{\boldsymbol{L}_{b}\boldsymbol{m}(\boldsymbol{\theta}\mathbf{)}%
=\boldsymbol{d}_{b}\}_{b=1}^{B}$\ so that the $(\boldsymbol{L}_{b}%
{,}\boldsymbol{d}_{b})$ is a submatrix of $(\boldsymbol{L}_{b+1}%
{,}\boldsymbol{d}_{b+1})$. Such constraints define a sequence of LMLC that
share the same design matrix $\boldsymbol{X}{=(}\boldsymbol{L}_{b}%
{,}\boldsymbol{W}_{b}{)}$.

\begin{Definition}
\label{Def7}The sequence of LMLC%
\begin{equation}
\log\boldsymbol{m}(\boldsymbol{\theta}\mathbf{)=}\boldsymbol{X}%
\boldsymbol{\theta}\quad\text{and}\quad\boldsymbol{L}_{b}^{T}\boldsymbol{m}%
(\boldsymbol{\theta}\mathbf{)}=\boldsymbol{d}_{b}{,}\label{eq28}%
\end{equation}
where $\boldsymbol{X}=(\boldsymbol{x}{_{1},...,}\boldsymbol{x}{_{t}})$,
$\boldsymbol{L}_{1}=(\boldsymbol{x}{_{1},...,}\boldsymbol{x}{_{c+r}})$ and
$\boldsymbol{L}_{b+1}=(\boldsymbol{L}_{b},\boldsymbol{x}{_{c+r+b}})$,
$b\in\{1,...,t-c-r\}$, is called the $b$-th \textsf{reduced LMLC through the
constraints}, because by increasing one unit the number of constraints, since
$\Theta_{b+1}\subset\Theta_{b}$ with the parameter space given by $\Theta
_{b}\equiv\{\boldsymbol{\theta}\in\mathbf{{%
\mathbb{R}
}}^{t}:\boldsymbol{L}_{b}^{T}\boldsymbol{m}(\boldsymbol{\theta}\mathbf{)=}%
\boldsymbol{d}_{b}\}$, it holds $M_{b+1}\subset M_{b}$\ where%
\[
M_{b}\equiv\{\boldsymbol{m}{(}\boldsymbol{\theta}\mathbf{)}\in\mathbf{{%
\mathbb{R}
}}^{k}:\log\boldsymbol{m}{(}\boldsymbol{\theta}\mathbf{)}=\boldsymbol{X}%
\boldsymbol{\theta},\boldsymbol{\theta}\in\Theta_{b}\}.
\]

\end{Definition}

In the following a generalized definition of nested LMLC is given, in which
Definitions \ref{Def6}\ and \ref{Def7} are covered.

\begin{Definition}
\label{Def8}In a sequence of LMLC $\{M_{b}\}_{b=1}^{B}$ such that%
\begin{align*}
M_{b} &  \equiv\{\boldsymbol{m}{(}\boldsymbol{\theta}_{b}\mathbf{)}%
\in\mathbf{{%
\mathbb{R}
}}^{k}:\log\boldsymbol{m}{(}\boldsymbol{\theta}_{b}\mathbf{)}=\boldsymbol{X}%
_{b}\boldsymbol{\theta}_{b},\boldsymbol{\theta}_{b}\in\Theta_{b}\},\\
\Theta_{b} &  \equiv\{\boldsymbol{\theta}_{b}\in\mathbf{{%
\mathbb{R}
}}^{t_{b}}:\boldsymbol{L}_{b}^{T}\boldsymbol{m}(\boldsymbol{\theta}%
_{b}\mathbf{)=}\boldsymbol{d}_{b}\},\\
t_{b} &  \equiv\operatorname*{rank}(\boldsymbol{X}_{b}),\\
\boldsymbol{L}_{b} &  \equiv(\boldsymbol{X}_{0},\boldsymbol{C}_{b}),\\
r_{b} &  \equiv\operatorname*{rank}(\boldsymbol{C}_{b}),
\end{align*}
$M_{b+1}$\ is said to be nested within $M_{b}$\ ($b\in\{1,...,B-1\}$), denoted
by $M_{b+1}\subset M_{b}$,if it holds%
\begin{equation}
\mathcal{C}(\boldsymbol{X}_{b+1})\subset\mathcal{C}(\boldsymbol{X}_{b}%
)\quad\text{and}\quad\mathcal{C}(\boldsymbol{L}_{b})\subset\mathcal{C}%
(\boldsymbol{L}_{b+1}),\label{eq31}%
\end{equation}
with $t_{b+1}\leq t_{b}$ and $r_{b+1}\geq r_{b}$, being strict at least one of
the two inequalities.
\end{Definition}

Once a sequence of nested LMLC $\{M_{b}\}_{b=1}^{B}$ has been established, our
goal is to present $\phi$-divergence test statistics to test successively
\begin{equation}
H_{Null}(b):M_{b+1}\text{ against }H_{Alt}(b):M_{b}-M_{b+1};\quad
b=1,...,B-1,\label{eq32}%
\end{equation}
where we continue to test as long as the null hypothesis is accepted and we
infer an integer $b^{\ast}$, such that $b\in\{1,...,B-1\}$, to be the first
value $b$ for which $M_{b+1}$\ is rejected as null hypothesis, or $b^{\ast}=B$ otherwise.

In Agresti \cite[Section 4.5.4]{Agresti} the classical likelihood ratio test
statistic for loglinear models ($r_{b}=r_{b+1}=0$, $c\geq0$) is given,%
\begin{align}
G^{2}(\widehat{\boldsymbol{\theta}}_{b+1}|\widehat{\boldsymbol{\theta}}_{b})
&  =2%
{\displaystyle\sum\limits_{{i}=1}^{k}}
\left(  {m}_{{i}}(\widehat{\boldsymbol{\theta}}_{b})\log\frac{{m}_{{i}%
}(\widehat{\boldsymbol{\theta}}_{b})}{{m}_{{i}}(\widehat{\boldsymbol{\theta}%
}_{b+1})}-{m}_{{i}}(\widehat{\boldsymbol{\theta}}_{b})+{m}_{{i}}%
(\widehat{\boldsymbol{\theta}}_{b+1})\right) \nonumber\\
&  =2\left(  D_{Kull}(\boldsymbol{n},\boldsymbol{m}(\widehat
{{\boldsymbol{\theta}}}_{b+1}))-D_{Kull}(\boldsymbol{n},\boldsymbol{m}%
(\widehat{{\boldsymbol{\theta}}}_{b}))\right)  ,\label{eq33}%
\end{align}
where $\widehat{\boldsymbol{\theta}}_{b+1}$ and $\widehat{\boldsymbol{\theta}%
}_{b}$ are the MLE's of the parameter\ in the models $M_{b+1}$\ and $M_{b}%
$\ respectively. It is also shown that $G^{2}(\widehat{\boldsymbol{\theta}%
}_{b+1}|\widehat{\boldsymbol{\theta}}_{b})$ is asymptotically distributed
according to a chi-square with $t_{b}-t_{b+1}$ degrees of freedom under the
null hypothesis of (\ref{eq32}). Minimizing the Kullback divergence measure
over a smaller parameter space cannot yield a larger minimum value, therefore
$D_{Kull}(\boldsymbol{n},\boldsymbol{m}(\widehat{{\boldsymbol{\theta}}}%
_{b+1}))\geq D_{Kull}(\boldsymbol{n},\boldsymbol{m}(\widehat
{{\boldsymbol{\theta}}}_{b}))$. However there is another interesting way to
show the same inequality, which is based on proving%
\begin{equation}
D_{Kull}(\boldsymbol{n},\boldsymbol{m}(\widehat{{\boldsymbol{\theta}}}%
_{b+1}))-D_{Kull}(\boldsymbol{n},\boldsymbol{m}(\widehat{{\boldsymbol{\theta}%
}}_{b}))=D_{Kull}\left(  \boldsymbol{m}(\widehat{{\boldsymbol{\theta}}}%
_{b}),\boldsymbol{m}(\widehat{{\boldsymbol{\theta}}}_{b+1})\right)
,\label{eq34}%
\end{equation}
whose non-negativity is guaranteed by the common property of all divergence
measures. Based on (\ref{eq33}) and (\ref{eq34}), in Mart\'{\i}n and Pardo
\cite[Section 4]{MartinPardo} divergence-based test statistics were introduced
for the same models ($r_{b}=r_{b+1}=0$, $c\geq0$),%
\begin{equation}
S^{\phi}(\widehat{\boldsymbol{\theta}}_{b+1}^{\phi}|\widehat
{\boldsymbol{\theta}}_{b}^{\phi})=\frac{2}{\phi^{\prime\prime}(1)}\left(
D_{\phi}(\boldsymbol{n},\boldsymbol{m}(\widehat{{\boldsymbol{\theta}}}%
_{b+1}^{\phi}))-D_{\phi}(\boldsymbol{n},\boldsymbol{m}(\widehat
{{\boldsymbol{\theta}}}_{b}^{\phi}))\right) \label{eq35}%
\end{equation}
and%
\begin{equation}
T^{\phi_{1}}(\widehat{\boldsymbol{\theta}}_{b+1}^{\phi_{2}}|\widehat
{\boldsymbol{\theta}}_{b}^{\phi_{2}})=\frac{2}{\phi_{1}^{\prime\prime}%
(1)}D_{\phi_{1}}\left(  \boldsymbol{m}(\widehat{{\boldsymbol{\theta}}}%
_{b}^{\phi_{2}}),\boldsymbol{m}(\widehat{{\boldsymbol{\theta}}}_{b+1}%
^{\phi_{2}})\right)  ,\label{eq36}%
\end{equation}
whose asymptotic distribution under the null hypothesis of (\ref{eq32}) was
shown to be exactly the same as $G^{2}(\widehat{\boldsymbol{\theta}}%
_{b+1}|\widehat{\boldsymbol{\theta}}_{b})$ for both of them. It should be
pointed out that (\ref{eq34}) does not hold by replacing any $\phi$-divergence
measure instead of the Kullback divergence measure.

In the more general framework of LMLC ($r_{b+1}\geq r_{b}\geq0$, $c\geq0$) we
shall establish herein that under the null hypothesis of (\ref{eq32}), the
test statistics $T^{\phi_{1}}(\widehat{\boldsymbol{\theta}}_{b+1}^{\phi_{2}%
}|\widehat{\boldsymbol{\theta}}_{b}^{\phi_{2}})$ and $S^{\phi}(\widehat
{\boldsymbol{\theta}}_{b+1}^{\phi}|\widehat{\boldsymbol{\theta}}_{b}^{\phi}%
)$\ converge in law to a chi-square with $t_{b}-t_{b+1}-r_{b}+r_{b+1}$ degrees
of freedom ($\chi_{t_{b}-t_{b+1}-r_{b}+r_{b+1}}^{2}$), $b=1,...,B-1$. Thus,
$\chi_{t_{b}-t_{b+1}-r_{b}+r_{b+1}}^{2}\left(  1-\alpha\right)  $ could be
chosen as a cutpoint for the rejection region.

\begin{theorem}
\label{Th9}Suppose that the data $\boldsymbol{n}=\left(  n_{1},...,n_{k}%
\right)  ^{T}$ are Poisson distributed. Choose the function $\phi_{1},\phi
_{2}\in\Phi$, where $\Phi$ was defined in Section \ref{Sec1}. Then, for
testing (\ref{eq32}), the asymptotic null distribution of the test statistic
$T^{\phi_{1}}(\widehat{\boldsymbol{\theta}}_{b+1}^{\phi_{2}}|\widehat
{\boldsymbol{\theta}}_{b}^{\phi_{2}})$, given in (\ref{eq36}), is chi-squared
with $t_{b}-t_{b+1}-r_{b}+r_{b+1}$ degrees of freedom.
\end{theorem}

\begin{proof}
A similar Taylor's expansion to one given in Theorem \ref{Th7} yields%
\[
T^{\phi_{1}}(\widehat{\boldsymbol{\theta}}_{b+1}^{\phi_{2}}|\widehat
{\boldsymbol{\theta}}_{b}^{\phi_{2}})=\boldsymbol{Z}_{b}^{T}\boldsymbol{Z}%
_{b}{+o}_{P}\left(  1\right)  ,
\]
where%
\[
\boldsymbol{Z}_{b}=\boldsymbol{D}_{\boldsymbol{m}{(}\boldsymbol{\theta}_{0}%
)}^{-\frac{1}{2}}(\boldsymbol{m}(\widehat{\boldsymbol{\theta}}_{b}^{\phi_{2}%
}){-}\boldsymbol{m}(\widehat{\boldsymbol{\theta}}_{b+1}^{\phi_{2}}))
\]
is distributed asymptotically as a normal distribution with mean vector zero
and variance-covariance matrix $\boldsymbol{T}_{b}^{\ast}=\boldsymbol{K}%
_{b}-\boldsymbol{K}_{b+1}$ with%
\begin{align}
\boldsymbol{K}_{j} &  {\equiv}\boldsymbol{A}_{X_{j}}{(\boldsymbol{\theta
}{_{b+1,0}})-}\boldsymbol{A}{_{X_{j}}{(\boldsymbol{\theta}{_{b+1,0}})}%
}\boldsymbol{D}_{\boldsymbol{m}^{\ast}(\boldsymbol{\theta}_{b+1,0}\mathbf{)}%
}^{\frac{1}{2}}\boldsymbol{L}_{j}\left(  \boldsymbol{L}_{j}^{T}\boldsymbol{D}%
_{\boldsymbol{m}^{\ast}(\boldsymbol{\theta}_{b+1,0}\mathbf{)}}^{\frac{1}{2}%
}\right. \label{eq39}\\
&  \times\left.  \boldsymbol{A}_{X_{j}}{(\boldsymbol{\theta}{_{b+1,0}}%
)}\boldsymbol{D}_{\boldsymbol{m}^{\ast}(\boldsymbol{\theta}_{b+1,0}\mathbf{)}%
}^{\frac{1}{2}}\boldsymbol{L}_{j}\right)  ^{-1}\boldsymbol{L}_{j}%
^{T}\boldsymbol{D}_{\boldsymbol{m}^{\ast}(\boldsymbol{\theta}_{b,0}\mathbf{)}%
}^{\frac{1}{2}}\boldsymbol{A}_{X_{j}}{(\boldsymbol{\theta}{_{b+1,0}}%
)},\nonumber
\end{align}%
\begin{equation}
\boldsymbol{A}_{X_{j}}{(\boldsymbol{\theta}{_{b+1,0}})}{\equiv}\boldsymbol{D}%
_{\boldsymbol{m}^{\ast}(\boldsymbol{\theta}_{b+1,0})}^{\frac{1}{2}%
}\boldsymbol{X}_{j}\!\left(  \boldsymbol{X}_{j}^{T}\!\boldsymbol{D}%
_{\boldsymbol{m}^{\ast}(\boldsymbol{\theta}_{b+1,0})}\boldsymbol{X}%
_{j}\right)  ^{-1}\!\boldsymbol{X}_{j}^{T}\!\boldsymbol{D}_{\boldsymbol{m}%
^{\ast}(\boldsymbol{\theta}_{b+1,0})}^{\frac{1}{2}},\label{eq38}%
\end{equation}
for $j=b,b+1$. The asymptotic distribution of the test statistic $T^{\phi_{1}%
}(\widehat{\boldsymbol{\theta}}_{b+1}^{\phi_{2}}|\widehat{\boldsymbol{\theta}%
}_{b}^{\phi_{2}})$ will be a chi-square if the matrix $\boldsymbol{T}%
_{b}^{\ast}$ is idempotent and symmetric. It is clear that $\boldsymbol{T}%
_{b}^{\ast}$ is symmetric, we shall establish that it is also idempotent.
Since $\mathcal{C(}\boldsymbol{D}_{\boldsymbol{m}^{\ast}(\boldsymbol{\theta
}_{b+1,0}\mathbf{)}}^{1/2}\boldsymbol{X}_{b+1}{)}\subset\mathcal{C(}%
\boldsymbol{D}_{\boldsymbol{m}^{\ast}(\boldsymbol{\theta}_{b+1,0}\mathbf{)}%
}^{1/2}\boldsymbol{X}_{b}{)}$\ we have%
\[
\boldsymbol{A}_{X_{b+1}}{(}\boldsymbol{\theta}{_{b+1,0})}{=}\boldsymbol{A}%
_{X_{b+1}}{(}\boldsymbol{\theta}{_{b+1,0})}\boldsymbol{A}_{X_{b+1}}%
{(}\boldsymbol{\theta}{_{b+1,0})=}\boldsymbol{A}_{X_{b+1}}{(}%
\boldsymbol{\theta}{_{b+1,0})}\boldsymbol{A}_{X_{b}}{(}\boldsymbol{\theta
}{_{b+1,0}),}%
\]%
\[
\boldsymbol{A}_{X_{b}}{(}\boldsymbol{\theta}{_{b+1,0})=}\boldsymbol{A}_{X_{b}%
}{(}\boldsymbol{\theta}{_{b+1,0})}\boldsymbol{A}_{X_{b}}{(}\boldsymbol{\theta
}{_{b+1,0}),}%
\]
and on the other hand since $\mathcal{C(}\boldsymbol{L}_{b}{)}\subset
\mathcal{C(}\boldsymbol{L}_{b+1}{)}$ there exists a matrix $\boldsymbol{B}%
$\ such that $\boldsymbol{L}_{b}=\boldsymbol{L}_{b+1}\boldsymbol{B}$. Thus it
holds\newline i) $\boldsymbol{K}_{b+1}=\boldsymbol{K}_{b+1}\boldsymbol{K}%
_{b+1}=\boldsymbol{K}_{b+1}\boldsymbol{K}_{b}=\boldsymbol{K}_{b}%
\boldsymbol{K}_{b+1} $,\newline ii) $\boldsymbol{K}_{b}=\boldsymbol{K}%
_{b}\boldsymbol{K}_{b}$,\newline which implies $\boldsymbol{T}_{b}^{\ast
}\boldsymbol{T}_{b}^{\ast}=\boldsymbol{T}_{b}^{\ast}$.\newline The degrees of
freedom of the chi-squared distributed statistic $T^{\phi_{1}}(\widehat
{\boldsymbol{\theta}}_{b+1}^{\phi_{2}}|\widehat{\boldsymbol{\theta}}_{b}%
^{\phi_{2}})$ coincides with the trace of the matrix $\boldsymbol{T}_{b}%
^{\ast}$, i.e. $t_{b}-t_{b+1}-r_{b}+r_{b+1}$. $\square$
\end{proof}

For the test statistic $S^{\phi}(\widehat{\boldsymbol{\theta}}_{b+1}^{\phi
}|\widehat{\boldsymbol{\theta}}_{b}^{\phi})$ the same result as Theorem
\ref{Th9}\ can be obtained by following a similar proof.

\begin{Remark}
\label{Rem4}Consider the \textsf{saturated LMLC}, i.e. the design matrix of
the loglinear model is given by a $k\times k$ matrix $\boldsymbol{X}_{1}%
$\ ($t_{1}=k$), and we may assume, without any loss of generality,
$\boldsymbol{X}_{1}{=}\boldsymbol{I}_{k}$. On the other hand, apart from the
constraints associated with the sampling scheme ($c\geq0$) there is no any
additional linear constraint ($r_{1}=0$)%
\[
M_{1}=\{\boldsymbol{m}{(}\boldsymbol{\theta}_{1}\mathbf{)}\in\mathbf{{%
\mathbb{R}
}}^{k}:\log\boldsymbol{m}{(}\boldsymbol{\theta}_{1}\mathbf{)}%
=\boldsymbol{\theta}_{1},\boldsymbol{\theta}_{1}\in\Theta_{1}\},
\]%
\[
\Theta_{1}=\{\boldsymbol{\theta}_{1}\in\mathbf{{%
\mathbb{R}
}}^{k}:\boldsymbol{X}_{0}^{T}\boldsymbol{m}(\boldsymbol{\theta}_{1}%
\mathbf{)=}\boldsymbol{X}_{0}^{T}\boldsymbol{n}\}.
\]
Consider also a generic LMLC%
\[
M_{2}=\{\boldsymbol{m}{(}\boldsymbol{\theta}_{2}\mathbf{)}\in\mathbf{{%
\mathbb{R}
}}^{k}:\log\boldsymbol{m}{(}\boldsymbol{\theta}_{2}\mathbf{)}=\boldsymbol{X}%
_{2}\boldsymbol{\theta}_{2},\boldsymbol{\theta}_{2}\in\Theta_{2}\},
\]%
\[
\Theta_{2}=\{\boldsymbol{\theta}_{2}\in\mathbf{{%
\mathbb{R}
}}^{t_{2}}:\boldsymbol{L}_{2}^{T}\boldsymbol{m}(\boldsymbol{\theta}%
_{2}\mathbf{)=}\boldsymbol{d}_{2}\},
\]
where $t_{2}\leq k$, $r_{2}\geq0$, being strict at least one of the two
inequalities. The hypothesis testing (\ref{eq32}) for the two$\ $nested LMLC
above ($B=2$) is the same as the goodness-of-fit test (\ref{eq24}) associated
with the model $M_{2}$. Therefore $T^{\phi_{1}}(\widehat{\boldsymbol{\theta}%
}_{2}^{\phi_{2}}|\widehat{\boldsymbol{\theta}}_{1}^{\phi_{2}})=T^{\phi_{1}%
}(\widehat{\boldsymbol{\theta}}_{2}^{\phi_{2}})$.
\end{Remark}

To test the sequence of LMLC (\ref{eq32}) $b=1,...,b^{\ast}$, we need an
asymptotic independence result for the sequence of test statistics
$\{T^{\phi_{1}}(\widehat{\boldsymbol{\theta}}_{b+1}^{\phi_{2}}|\widehat
{\boldsymbol{\theta}}_{b}^{\phi_{2}})\}_{b=1}^{b^{\ast}}$ (or $\{S^{\phi
}(\widehat{\boldsymbol{\theta}}_{b+1}^{\phi}|\widehat{\boldsymbol{\theta}}%
_{b}^{\phi})\}_{b=1}^{b^{\ast}}$). This result is given in the theorem below.

\begin{theorem}
\label{Th10}Suppose that data $\boldsymbol{n}=\left(  n_{1},...,n_{k}\right)
^{T}$ are Poisson distributed. We first test, $H_{Null}:M_{b+1}$ against
$H_{Alt}:M_{b},$ followed by $H_{Null}:M_{b}$ against $H_{Alt}:M_{b-1}$. Then,
under the assumption that it holds $M_{b+1}$, the statistics $T^{\phi_{1}%
}(\widehat{\boldsymbol{\theta}}_{b+1}^{\phi_{2}}|\widehat{\boldsymbol{\theta}%
}_{b}^{\phi_{2}})$ and $T^{\phi_{1}}(\widehat{\boldsymbol{\theta}}_{b}%
^{\phi_{2}}|\widehat{\boldsymbol{\theta}}_{b-1}^{\phi_{2}})$\ are
asymptotically independent.
\end{theorem}

\begin{proof}
A second order Taylor's expansion gives%
\[
T^{\phi_{1}}(\widehat{\boldsymbol{\theta}}_{j}^{\phi_{2}}|\widehat
{\boldsymbol{\theta}}_{j-1}^{\phi_{2}})=\widetilde{\boldsymbol{Z}}_{j}%
^{T}\widetilde{\boldsymbol{Z}}_{j}{+}o_{P}(1),\qquad j\in\{b+1,b\},
\]
where%
\[
\widetilde{\boldsymbol{Z}}_{j}=\boldsymbol{T}_{j-1}^{\ast}\boldsymbol{D}%
_{\boldsymbol{m}^{\ast}(\boldsymbol{\theta}_{b+1,0}\mathbf{)}}^{-\frac{1}{2}%
}{\frac{1}{\sqrt{N}}(}\boldsymbol{n}{-}\boldsymbol{m}{({\boldsymbol{\theta}%
}_{b+1,0})),}%
\]
and $\boldsymbol{T}_{j}^{\ast}=\boldsymbol{K}_{j}-\boldsymbol{K}_{j+1}$ with
$\boldsymbol{K}_{j}$,\ $j=b,b-1$ defined in (\ref{eq39}). By Searle
\cite[Theorem 4 in page 59]{Searle} the quadratic forms%
\[
\widetilde{\boldsymbol{Z}}_{j}^{T}\widetilde{\boldsymbol{Z}}_{j}={\tfrac
{1}{\sqrt{N}}(}\boldsymbol{n}{-}\boldsymbol{m}{({\boldsymbol{\theta}}%
_{b+1,0}))}^{T}\boldsymbol{D}_{\boldsymbol{m}^{\ast}(\boldsymbol{\theta
}_{b+1,0}\mathbf{)}}^{-\frac{1}{2}}\boldsymbol{T}_{j-1}^{\ast}\boldsymbol{D}%
_{\boldsymbol{m}^{\ast}(\boldsymbol{\theta}_{b+1,0}\mathbf{)}}^{-\frac{1}{2}%
}{\tfrac{1}{\sqrt{N}}(}\boldsymbol{n}{-}\boldsymbol{m}{({\boldsymbol{\theta}%
}_{b+1,0}))},
\]
for $j=b+1,b$, are asymptotically independent if
\[
\boldsymbol{D}_{\boldsymbol{m}^{\ast}(\boldsymbol{\theta}_{b+1,0}\mathbf{)}%
}^{-\frac{1}{2}}\boldsymbol{T}_{b+1}^{\ast}\boldsymbol{D}_{\boldsymbol{m}%
^{\ast}(\boldsymbol{\theta}_{b+1,0}\mathbf{)}}^{-\frac{1}{2}}%
\boldsymbol{\Sigma}\boldsymbol{D}_{\boldsymbol{m}^{\ast}(\boldsymbol{\theta
}_{b+1,0}\mathbf{)}}^{-\frac{1}{2}}\boldsymbol{T}_{b}^{\ast}\boldsymbol{D}%
_{\boldsymbol{m}^{\ast}(\boldsymbol{\theta}_{b+1,0}\mathbf{)}}^{-\frac{1}{2}%
}=\boldsymbol{0}_{k\times k}{,}%
\]
where matrix $\boldsymbol{\Sigma}=\boldsymbol{D}_{\boldsymbol{m}^{\ast
}(\boldsymbol{\theta}_{b+1,0}\mathbf{)}}$\ is the asymptotic
variance-covariance matrix of vector ${(}\boldsymbol{n}{-}\boldsymbol{m}{(}%
${{$\boldsymbol{\theta}$}}${_{b+1,0}))/}\sqrt{N}$. By following a similar
argument given in the proof of Theorem \ref{Th9} to see that $\boldsymbol{T}%
_{b}^{\ast}$ is idempotent, it follows that $\boldsymbol{T}_{b+1}^{\ast
}\boldsymbol{T}_{b}^{\ast}=\boldsymbol{0}_{k\times k}$. $\square$
\end{proof}

For (\ref{eq35}) the same result as Theorem \ref{Th9}\ can be obtained by
following a similar proof.

According to Agresti \cite[page 215]{Agresti} once an asymptotic probability
of error I equals $1-(1-\alpha)^{\frac{1}{B-1}}$ has been established for each
test in a sequence of nested tests, the overall asymptotic probability of type
I error is less or equal than $\alpha$. In the next theorem a stronger result
is given.

\begin{theorem}
\label{Th12}For a sequence of $B-1$ tests (\ref{eq32}) associated with a
sequence of LMLC $\{M_{b}\}_{b=1}^{B}$, when each test has a size equals
$1-(1-\alpha)^{\frac{1}{B-1}}$, the overall size of the tests is given by
$\alpha$.
\end{theorem}

\begin{proof}
For the purpose of establishing a size equals $1-(1-\alpha)^{\frac{1}{B-1}}$
for each hypothesis testing in (\ref{eq27}) we shall consider according to
Theorem \ref{Th9} a cutpoint for the rejection region equals $\chi_{df}%
^{2}((1-\alpha)^{\frac{1}{B-1}})$, where $df=t_{b}-t_{b+1}-r_{b}+r_{b+1}$.
Thus, the overall size for testing (\ref{eq32})\ in a sequence of nested LMLC
$\{M_{b}\}_{b=1}^{B}$ is given by%
\begin{align*}
& \Pr(\exists b\in\{1,...,B-1\}:T^{\phi_{1}}(\widehat{\boldsymbol{\theta}%
}_{b+1}^{\phi_{2}}|\widehat{\boldsymbol{\theta}}_{b}^{\phi_{2}})>\chi_{df}%
^{2}((1-\alpha)^{\frac{1}{B-1}})|H_{Null}(b))\\
& =1-\Pr(T^{\phi_{1}}(\widehat{\boldsymbol{\theta}}_{b+1}^{\phi_{2}}%
|\widehat{\boldsymbol{\theta}}_{b}^{\phi_{2}})\leq\chi_{df}^{2}((1-\alpha
)^{\frac{1}{B-1}})|H_{Null}(b),b=1,...,B-1)\\
& =1-%
{\textstyle\prod_{b=1}^{B-1}}
\Pr(T^{\phi_{1}}(\widehat{\boldsymbol{\theta}}_{b+1}^{\phi_{2}}|\widehat
{\boldsymbol{\theta}}_{b}^{\phi_{2}})\leq\chi_{df}^{2}((1-\alpha)^{\frac
{1}{B-1}})|H_{Null}(b))\\
& =1-%
{\textstyle\prod_{b=1}^{B-1}}
\Pr(\chi_{df}^{2}\leq\chi_{df}^{2}((1-\alpha)^{\frac{1}{B-1}})))\\
& =1-%
{\textstyle\prod_{b=1}^{B-1}}
(1-\alpha)^{\frac{1}{B-1}}=\alpha.
\end{align*}
\newline The second equality comes from Theorem \ref{Th10} and the third one
from Theorem \ref{Th9}. $\square$
\end{proof}

\section{Simulation study: Marginal homogeneity\label{Sec5}}

\subsection{Description of conditional and unconditional tests}

The traditionally so-called conditional test for marginal homogeneity (MH) was
applied for the first time in Caussinus \cite{Caussinus}. He noted that once
it is known that the quasi-symmetry (QS) model holds, marginal homogeneity
(MH) is equivalent to symmetry (S). In other words, because QS is a nested
model within S ($M_{\mathcal{QS}}\subset M_{\mathcal{S}}$), first we could
test whether it holds QS model against the alternative hypothesis of saturated
model (SAT), defined in Remark \ref{Rem4} ($M_{\mathcal{QS}}\subset
M_{\mathcal{SAT}}$),%
\begin{equation}
H_{Null}(1):M_{\mathcal{QS}}\text{ against }H_{Alt}(1):M_{\mathcal{SAT}%
}-M_{\mathcal{QS}},\label{eq30}%
\end{equation}
and after that%
\begin{equation}
H_{Null}(2):M_{\mathcal{S}}\text{ against }H_{Alt}(2):M_{\mathcal{S}%
}-M_{\mathcal{QS}}.\label{eq42}%
\end{equation}
Focussed on a square $I\times I$\ contingency table with multinomial sampling
($c=1$), one could be interested in analyzing what the difference is between
testing the conditional model above and the unconditional model of MH below%
\begin{equation}
H_{Null}:M_{\mathcal{MH}}\text{ against }H_{Alt}:M_{\mathcal{SAT}%
}-M_{\mathcal{MH}}.\label{eq43}%
\end{equation}
The formulation of these models for two-way contingency tables is%
\[%
{\displaystyle\sum\limits_{j=1}^{I}}
m_{ij}({\boldsymbol{\theta}}_{\mathcal{MH}})=%
{\displaystyle\sum\limits_{i=1}^{I}}
m_{ij}({\boldsymbol{\theta}}_{\mathcal{MH}})\text{ or\ }m_{i\bullet
}({\boldsymbol{\theta}}_{\mathcal{MH}})=m_{\bullet i}({\boldsymbol{\theta}%
}_{\mathcal{MH}})\text{,}\ i,j=1,...,I;
\]%
\[
\log m_{ij}({\boldsymbol{\theta}}_{\mathcal{S}})=u+\theta_{i}+\theta
_{j}+\theta_{ij},\quad i,j=1,...,I,
\]
with $\theta_{ij}=\theta_{ji}$, $\forall i\neq j$, $%
{\textstyle\sum\nolimits_{i=1}^{I}}
\theta_{i}=0$, $%
{\textstyle\sum\nolimits_{i=1}^{I}}
\theta_{12(ij)}=0$, $j=1,...,I$;%
\[
\log m_{ij}({\boldsymbol{\theta}_{\mathcal{OQS}}})=u+\theta_{1(i)}%
+\theta_{1(j)}+\beta w_{j}+\theta_{12(ij)},\quad i,j=1,...,I,
\]
with $\theta_{12(ij)}=\theta_{12(ji)}$, $\forall i\neq j$, $%
{\textstyle\sum\nolimits_{i=1}^{I}}
\theta_{1(i)}=0$, $%
{\textstyle\sum\nolimits_{i=1}^{I}}
\theta_{12(ij)}=0$, $j=1,...,I$, $%
{\textstyle\sum\nolimits_{j=1}^{I}}
w_{j}=0$, $%
{\textstyle\sum\nolimits_{j=1}^{I}}
w_{j}^{2}=1$, where $\{w_{j}\}_{j=1}^{I}$ is a set of weights associated to
each category $j\in\{1,...,I\}$ such that the distance between the contiguous
ones is fixed, i.e.
\[
w_{j}=\frac{j-\sum\limits_{i=1}^{I}i\left/  \overset{}{I}\right.  }{\left(
\sum\limits_{i=1}^{I}i^{2}-\frac{1}{I}\left(  \sum\limits_{i=1}^{I}i\right)
^{2}\right)  ^{1/2}}=\frac{2j-\left(  I+1\right)  }{\sqrt{\frac{1}{3}I\left(
I-1\right)  \left(  I+1\right)  }}.
\]
(for more details about the interpretation of this model see Agresti and
Kateri \cite{AgrestiKateri});%
\[
\log m_{ij}({\boldsymbol{\theta}_{\mathcal{QS}}})=u+\theta_{1(i)}%
+\theta_{2(j)}+\theta_{12(ij)},\quad i,j=1,...,I,
\]
with $\theta_{12(ij)}=\theta_{12(ji)}$, $\forall i\neq j$, $%
{\textstyle\sum\nolimits_{i=1}^{I}}
\theta_{1(i)}=%
{\textstyle\sum\nolimits_{j=1}^{I}}
\theta_{2(j)}=0$, $%
{\textstyle\sum\nolimits_{i=1}^{I}}
\theta_{12(ij)}=0$, $j=1,...,I$.

While $M_{\mathcal{S}}$ and$\ M_{\mathcal{QS}}$\ are loglinear models,
$M_{\mathcal{MH}}$ is a marginal model (see Remark \ref{Rem1}). By taking into
account the meaning of both tests, the initial conditions are different,
actually while a rejection of $H_{Null}$\ implies that MH is not accepted, a
rejection of $H_{Null}(1)$ does not implies the same fact, i.e. even though QS
is rejected a MH could hold. However, an acceptance of $H_{Null}$\ and
$H_{Null}(2)$\ implies the same hypothesis, i.e. MH is accepted. Because all
models, $M_{\mathcal{SAT}}$, $M_{\mathcal{QS}}$, $M_{\mathcal{S}}$, and
$M_{\mathcal{MH}}$, are LMLC, it is possible to establish the same true model
for both tests (\ref{eq30})-(\ref{eq42})\ and (\ref{eq43}), by choosing a
suitable parametrization according to the design matrices. On the other hand,
it will be possible to carry out such a test in order to compare the exact
size and power, through a simulation experiment by using $\phi$-divergence
based test-statistics. For this purpose we shall focus on M$\phi_{(\lambda
_{2})}$E's with $\lambda_{2}\in\{-0.5,0,2/3,1,2\}$ (i.e., with $\lambda_{2}%
=0$\ MLE's are included), as well as on the same family of $\phi_{(\lambda
_{1})}$-divergence measures (\ref{eq5}) for building test-statistics, with
$\lambda_{1}\in\{-0.5,0,2/3,1,2\}$ (i.e., once MLE's are included, with
$\lambda_{1}=0$\ and $\lambda_{1}=1$ the classical test-statistics are
obtained, likelihood ratio $G^{2}(\widehat{\boldsymbol{\theta}}_{b+1}%
|\widehat{\boldsymbol{\theta}}_{b})$ and chi-squared $X^{2}(\widehat
{\boldsymbol{\theta}}_{b+1}|\widehat{\boldsymbol{\theta}}_{b})$
test-statistics respectively)%
\[
T^{\phi_{(\lambda_{1})}}(\widehat{\boldsymbol{\theta}}_{b+1}^{\phi
_{(\lambda_{2})}}|\widehat{\boldsymbol{\theta}}_{b}^{\phi_{(\lambda_{2})}%
})=\left\{
\begin{array}
[c]{l}%
\frac{2}{\lambda_{1}(\lambda_{1}+1)}\!\left(
{\textstyle\sum\limits_{i=1}^{I}}
{\textstyle\sum\limits_{j=1}^{I}}
\frac{m_{ij}^{\lambda_{1}+1}(\widehat{\boldsymbol{\theta}}_{b-1}%
^{\phi_{(\lambda_{2})}})}{m_{ij}^{\lambda_{1}}(\widehat{\boldsymbol{\theta}%
}_{b}^{\phi_{(\lambda_{2})}})}-n\right)  \!,\lambda_{1}(\lambda_{1}+1)\neq0\\
2%
{\textstyle\sum\limits_{i=1}^{I}}
{\textstyle\sum\limits_{j=1}^{I}}
m_{ij}(\widehat{\boldsymbol{\theta}}_{b-1}^{\phi_{(\lambda_{2})}})\log
\frac{m_{ij}(\widehat{\boldsymbol{\theta}}_{b-1}^{\phi_{(\lambda_{2})}}%
)}{m_{ij}(\widehat{\boldsymbol{\theta}}_{b}^{\phi_{(\lambda_{2})}})}%
,\quad\lambda_{1}=0\\
2%
{\textstyle\sum\limits_{i=1}^{I}}
{\textstyle\sum\limits_{j=1}^{I}}
m_{ij}(\widehat{\boldsymbol{\theta}}_{b}^{\phi_{(\lambda_{2})}})\log
\frac{m_{ij}(\widehat{\boldsymbol{\theta}}_{b}^{\phi_{(\lambda_{2})}})}%
{m_{ij}(\widehat{\boldsymbol{\theta}}_{b-1}^{\phi_{(\lambda_{2})}})}%
,\quad\lambda_{1}=-1
\end{array}
\right.
\]
where $n$ and $I$ are respectively the total table count and the table size
(i.e., $k=I^{2}$ and $n=%
{\textstyle\sum\nolimits_{i=1}^{I}}
{\textstyle\sum\nolimits_{j=1}^{I}}
n_{ij}$). In particular, if the saturated LMLC is considered as $M_{b}$, then
$m_{ij}(\widehat{\boldsymbol{\theta}}_{b}^{\phi_{(\lambda_{2})}})=n_{ij}$ (see
Remark \ref{Rem4}).

We shall also consider another conditional test for MH, which is based on the
ordinal quasi-symmetry model (OQS), instead of the previously considered QS,%
\begin{equation}
H_{Null}(1^{\prime}):M_{\mathcal{OQS}}\text{ against }H_{Alt}(1^{\prime
}):M_{\mathcal{SAT}}-M_{\mathcal{OQS}},\label{eq44}%
\end{equation}
and%
\begin{equation}
H_{Null}(2^{\prime}):M_{\mathcal{S}}\text{ against }H_{Alt}(2^{\prime
}):M_{\mathcal{S}}-M_{\mathcal{OQS}}\label{eq45}%
\end{equation}
($M_{\mathcal{S}}\subset M_{\mathcal{OQS}}\subset M_{\mathcal{SAT}}$).
Although OQS is usually applied for ordinal categorical data, because its
interpretation, it is possible to consider OQS as LMLC, in a generic way, by
defining its design matrix (for more information about this model see Agresti
\cite[Section 8.4]{Agresti2}). In order to compute the powers for the two
conditional tests, apart from considering a common true model, because
$M_{\mathcal{S}}\subset M_{\mathcal{OQS}}\subset M_{\mathcal{QS}}$, we shall
consider the same points of the alternative hypotheses.

\begin{table}[bh]
\caption{Theoretical probabilities for a $I\times I$ table ($I=4$)}%
\label{Table1}
\begin{center}%
\begin{tabular}
[c]{c|cccc|c}%
$m_{ij}^{\ast}(\hbox{\boldmath$\theta$}_{0})$ & $1$ & $2$ & $3$ & $4$ &
$m_{i\bullet}^{\ast}(\hbox{\boldmath$\theta$}_{0})$\\\hline
$1$ & $0.08161$ & $0.03156$ & $0.01647$ & $0.01050$ & $0.14017$\\
$2$ & $0.03156$ & $0.21104$ & $0.05204$ & $0.01418$ & $0.30883$\\
$3$ & $0.01647$ & $0.05204$ & $0.22186$ & $0.03156$ & $0.32195$\\
$4$ & $0.01050$ & $0.01418$ & $0.03156$ & $0.17278$ & $0.22905$\\\hline
$m_{\bullet j}^{\ast}(\hbox{\boldmath$\theta$}_{0})$ & $0.14017$ & $0.30883$ &
$0.32195$ & $0.22905$ & $m_{\bullet\bullet}^{\ast}(\hbox{\boldmath$\theta
$}_{0})=1$%
\end{tabular}
\end{center}
\end{table}

In Table \ref{Table1} the theoretical probability vector belonging to a
multinomial sampling scheme with $n\in\{100,250,400,550\}$ is shown. Its
corresponding values for the parameters for each model (null hypotheses) are
also given ($t_{\mathcal{MH}}=16$, $t_{\mathcal{S}}=10$, $t_{\mathcal{OQS}%
}=11$ and $t_{\mathcal{QS}}=13$):
\begin{align*}
{\boldsymbol{\theta}_{\mathcal{MH}}} &  =(u_{\mathcal{MH}},{\theta
_{1(1)},\theta_{1(2)},\theta_{1(3)},\theta_{2(1)},\theta_{2(2)},\theta
_{2(3)},\theta_{12(11)}},{\theta_{12(12)},\theta_{12(13)},}\\
&  \hspace*{0.5cm}{\theta_{12(21)}},{\theta_{12(22)},\theta_{12(23)}%
,\theta_{12(31)}},{\theta_{12(32)},\theta_{12(33)})}^{T}\\
&  =(u_{\mathcal{MH}},-0.95,-1.6,-2.05,-0.95,0.95,-0.45,-1.75,-1.6,-0.45,\\
&  \hspace*{0.5cm}1.0,-0.95,-2.05,-1.75,-0.95,0.75)^{T},
\end{align*}%
\begin{align*}
{\boldsymbol{\theta}_{\mathcal{S}}} &  =(u_{\mathcal{S}},{\theta_{1}%
,\theta_{2},\theta_{3},\theta_{11}},{\theta_{22},\theta_{12},\theta
_{13},\theta_{24},\theta_{34}})^{T}\\
&  =(u_{\mathcal{S}},-0.35,0.25,0.3,1.5,1.25,-0.05,-0.75,-1,-0.25)^{T},
\end{align*}%
\[
{\boldsymbol{\theta}_{\mathcal{OQS}}}=({\boldsymbol{\theta}_{\mathcal{S}}%
},0)^{T},{\boldsymbol{\theta}_{\mathcal{QS}}}=({\boldsymbol{\theta
}_{\mathcal{S}}},0,0,0)^{T}%
\]
(the values of $u_{\mathcal{MH}}$\ and $u_{\mathcal{S}}$, are obtained through
the different sampling sizes).

\subsection{Simulated Exact sizes and powers}

Focussing first on the tests (\ref{eq44})-(\ref{eq45}), the simulation study
is based on repeating the random experiments described above $R=10,000$ times
to compute on one hand the \textquotedblleft exact sizes\textquotedblright\ by
simulation%
\begin{align*}
\alpha_{n}^{(\lambda_{1},\lambda_{2})}(\boldsymbol{\theta}_{\mathcal{OQS}}) &
=\tfrac{\#\{T^{\phi_{(\lambda_{1})}}(\widehat{\boldsymbol{\theta}%
}_{\mathcal{OQS}}^{\phi_{(\lambda_{2})}}|\widehat{\boldsymbol{\theta}%
}_{\mathcal{SAT}}^{\phi_{(\lambda_{2})}})>\chi_{\frac{(I+1)(I-2)}{2}}%
^{2}((1-\alpha)^{\frac{1}{2}})|\boldsymbol{\theta}_{\mathcal{OQS}})\}}{R},\\
\alpha_{n}^{(\lambda_{1},\lambda_{2})}(\boldsymbol{\theta}_{\mathcal{S}}) &
=\tfrac{\#\{T^{\phi_{(\lambda_{1})}}(\widehat{\boldsymbol{\theta}%
}_{\mathcal{S}}^{\phi_{(\lambda_{2})}}|\widehat{\boldsymbol{\theta}%
}_{\mathcal{OQS}}^{\phi_{(\lambda_{2})}})>\chi_{1}^{2}((1-\alpha)^{\frac{1}%
{2}})|\boldsymbol{\theta}_{\mathcal{S}})\}}{R},\\
\alpha_{n}^{(\lambda_{1},\lambda_{2})} &  =1-(1-\alpha_{n}^{(\lambda
_{1},\lambda_{2})}(\boldsymbol{\theta}_{\mathcal{OQS}}))(1-\alpha
_{n}^{(\lambda_{1},\lambda_{2})}(\boldsymbol{\theta}_{\mathcal{S}})),
\end{align*}
once a simulated \textquotedblleft nominal\textquotedblright\ size of
$1-(1-\alpha)^{\frac{1}{2}}$, with $\alpha=0.05$, has been chosen for each
test. In order to calculate simulated exact powers $12$ points are chosen ($6
$ for (\ref{eq44}) and $6$ for (\ref{eq45}))%
\begin{align}
{\boldsymbol{\theta}_{\mathcal{SAT}}(i)} &  =(u_{\mathcal{SAT}},{\theta
_{1(1)},\theta_{1(2)},\theta_{1(3)},\theta_{2(1)},\theta_{2(2)},\theta
_{2(3)}{+}\delta_{1}(i),\theta_{12(11)}},{\theta_{12(12)},}\nonumber\\
&  \hspace*{0.5cm}{\theta_{12(13)}{+{\delta}}_{2}(i),\theta_{12(21)}}%
,{\theta_{12(22)},\theta_{12(23)},\theta_{12(31)}},{\theta_{12(32)}%
,\theta_{12(33)})}^{T}\nonumber\\
&  =(u_{\mathcal{SAT}},-0.95,-1.6,-2.05,-0.95,0.95,-0.45{+\delta}%
_{1}(i),-1.75,-1.6,\nonumber\\
&  \hspace*{0.5cm}-0.45{{+{\delta}}_{2}(i)}%
,1.0,-0.95,-2.05,-1.75,-0.95,0.75)^{T},\label{eq46}%
\end{align}
with $(({\delta}_{1}(1),{\delta}_{2}(1)),...,({\delta}_{1}(6),{\delta}%
_{2}(6)))^{T}=((0.45,0),(0.7,0),(0.9,0),(0,0.45),(0,0.7),(0,0.9))^{T}$.%
\begin{align}
{\boldsymbol{\theta}_{\mathcal{OQS}}(i)} &  =(u_{\mathcal{OQS}},{\theta
_{1},\theta_{2},\theta_{3},\theta_{11}},{\theta_{22},\theta_{12},\theta
_{13},\theta_{24},\theta_{34},\beta(i)})^{T}\nonumber\\
&  =(u_{\mathcal{OQS}},-0.35,0.25,0.3,1.5,1.25,-0.05,-0.75,-1,-0.25,{\beta
(i)})^{T},\label{eq47}%
\end{align}
with $({\beta(7),...,\beta(12)})^{T}=(0.5,0.7,1.0,-0.5,-0.7,-1.0)^{T}$ (the
values of $u_{\mathcal{SAT}}$ and $u_{\mathcal{OQS}}$\ are obtained through
the different sampling sizes). Thus the simulated exact powers are given by%
\[
\beta_{n}^{(\lambda_{1},\lambda_{2})}(\boldsymbol{\theta}_{\mathcal{SAT}%
}(i))=\tfrac{\#\{T^{\phi_{(\lambda_{1})}}(\widehat{\boldsymbol{\theta}%
}_{\mathcal{OQS}}^{\phi_{(\lambda_{2})}}|\widehat{\boldsymbol{\theta}%
}_{\mathcal{SAT}}^{\phi_{(\lambda_{2})}})>\chi_{\frac{(I+1)(I-2)}{2}}%
^{2}((1-\alpha)^{\frac{1}{2}})|\boldsymbol{\theta}_{\mathcal{SAT}}(i))\}}{R},
\]
$i=1,...,6$,
\[
\beta_{n}^{(\lambda_{1},\lambda_{2})}(\boldsymbol{\theta}_{\mathcal{OQS}%
}(i))=\tfrac{\#\{T^{\phi_{(\lambda_{1})}}(\widehat{\boldsymbol{\theta}%
}_{\mathcal{S}}^{\phi_{(\lambda_{2})}}|\widehat{\boldsymbol{\theta}%
}_{\mathcal{OQS}}^{\phi_{(\lambda_{2})}})>\chi_{1}^{2}((1-\alpha)^{\frac{1}%
{2}})|\boldsymbol{\theta}_{\mathcal{OQS}}(i))\}}{R},
\]
$i=7,...,12$. Focussing on the tests (\ref{eq30})-(\ref{eq42}), simulated
exact sizes are given by%
\begin{align*}
\alpha_{n}^{(\lambda_{1},\lambda_{2})}(\boldsymbol{\theta}_{\mathcal{QS}}) &
=\tfrac{\#\{T^{\phi_{(\lambda_{1})}}(\widehat{\boldsymbol{\theta}%
}_{\mathcal{QS}}^{\phi_{(\lambda_{2})}}|\widehat{\boldsymbol{\theta}%
}_{\mathcal{SAT}}^{\phi_{(\lambda_{2})}})>\chi_{\frac{(I-1)(I-2)}{2}}%
^{2}((1-\alpha)^{\frac{1}{2}})|\boldsymbol{\theta}_{\mathcal{QS}})\}}{R},\\
\alpha_{n}^{(\lambda_{1},\lambda_{2})}(\boldsymbol{\theta}_{\mathcal{S}}) &
=\tfrac{\#\{T^{\phi_{(\lambda_{1})}}(\widehat{\boldsymbol{\theta}%
}_{\mathcal{S}}^{\phi_{(\lambda_{2})}}|\widehat{\boldsymbol{\theta}%
}_{\mathcal{QS}}^{\phi_{(\lambda_{2})}})>\chi_{I-1}^{2}((1-\alpha)^{\frac
{1}{2}})|\boldsymbol{\theta}_{\mathcal{S}})\}}{R},\\
\alpha_{n}^{(\lambda_{1},\lambda_{2})} &  =1-(1-\alpha_{n}^{(\lambda
_{1},\lambda_{2})}(\boldsymbol{\theta}_{\mathcal{QS}}))(1-\alpha_{n}%
^{(\lambda_{1},\lambda_{2})}(\boldsymbol{\theta}_{\mathcal{S}})),
\end{align*}
and to calculate simulated exact powers $12$ points are chosen ($6$ for
(\ref{eq30}) and $6$ for (\ref{eq42})), the same points (\ref{eq46}) and
(\ref{eq47}) are valid taking into account that $\boldsymbol{\theta
}_{\mathcal{QS}}(i)=(\boldsymbol{\theta}_{\mathcal{OQS}}(i),0,0)^{T}$,%
\[
\beta_{n}^{(\lambda_{1},\lambda_{2})}(\boldsymbol{\theta}_{\mathcal{SAT}%
}(i))=\tfrac{\#\{T^{\phi_{(\lambda_{1})}}(\widehat{\boldsymbol{\theta}%
}_{\mathcal{QS}}^{\phi_{(\lambda_{2})}}|\widehat{\boldsymbol{\theta}%
}_{\mathcal{SAT}}^{\phi_{(\lambda_{2})}})>\chi_{\frac{(I-1)(I-2)}{2}}%
^{2}((1-\alpha)^{\frac{1}{2}})|\boldsymbol{\theta}_{\mathcal{SAT}}(i))\}}{R},
\]
$i=1,...,6$,
\[
\beta_{n}^{(\lambda_{1},\lambda_{2})}(\boldsymbol{\theta}_{\mathcal{QS}%
}(i))=\tfrac{\#\{T^{\phi_{(\lambda_{1})}}(\widehat{\boldsymbol{\theta}%
}_{\mathcal{S}}^{\phi_{(\lambda_{2})}}|\widehat{\boldsymbol{\theta}%
}_{\mathcal{QS}}^{\phi_{(\lambda_{2})}})>\chi_{I-1}^{2}((1-\alpha)^{\frac
{1}{2}})|\boldsymbol{\theta}_{\mathcal{QS}}(i))\}}{R},
\]
$i=7,...,12$. For the goodness of fit test (\ref{eq30}) we have, as usual,%
\[
\alpha_{n}^{(\lambda_{1},\lambda_{2})}(\boldsymbol{\theta}_{\mathcal{MH}%
})=\tfrac{\#\{T^{\phi_{(\lambda_{1})}}(\widehat{\boldsymbol{\theta}%
}_{\mathcal{MH}}^{\phi_{(\lambda_{2})}})>\chi_{I-1}^{2}((1-\alpha)^{\frac
{1}{2}})|\boldsymbol{\theta}_{\mathcal{MH}})\}}{R}%
\]
and to calculate simulated exact powers the same $12$ points above are chosen%
\[
\beta_{n}^{(\lambda_{1},\lambda_{2})}(\boldsymbol{\theta}_{\mathcal{SAT}%
}(i))=\tfrac{\#\{T^{\phi_{(\lambda_{1})}}(\widehat{\boldsymbol{\theta}%
}_{\mathcal{MH}}^{\phi_{(\lambda_{2})}})>\chi_{I-1}^{2}((1-\alpha)^{\frac
{1}{2}})|\boldsymbol{\theta}_{\mathcal{SAT}}(i))\}}{R},
\]
$i=1,...,6$,%
\[
\beta_{n}^{(\lambda_{1},\lambda_{2})}(\boldsymbol{\theta}_{\mathcal{OQS}%
}(i))=\tfrac{\#\{T^{\phi_{(\lambda_{1})}}(\widehat{\boldsymbol{\theta}%
}_{\mathcal{MH}}^{\phi_{(\lambda_{2})}}))>\chi_{I-1}^{2}((1-\alpha)^{\frac
{1}{2}})|\boldsymbol{\theta}_{\mathcal{OQS}}(i))\}}{R},
\]
$i=7,...,12$. The results of simulated exact sizes\ are shown in Tables
\ref{Table2} and \ref{Table3}. \begin{table}[ptbh]
\caption{$\alpha_{n}^{(\lambda_{1},\lambda_{2})}$ with $n\in\{100,250\}$}%
\label{Table2}
\begin{center}
$%
\begin{tabular}
[c]{cc|ccc|c|c|c|}\cline{3-8}
&  & \multicolumn{3}{|c|}{$n=100$} & \multicolumn{3}{|c|}{$n=250$}\\\hline
\multicolumn{1}{|c|}{$\lambda_{2}$} & \multicolumn{1}{|c|}{$\lambda_{1}$} &
\multicolumn{1}{|c|}{(\ref{eq43})} & \multicolumn{1}{c|}{(\ref{eq30}%
)-(\ref{eq42})} & (\ref{eq44})-(\ref{eq45}) & (\ref{eq43}) & (\ref{eq30}%
)-(\ref{eq42}) & (\ref{eq44})-(\ref{eq45})\\\hline
\multicolumn{1}{|c|}{} & \multicolumn{1}{|c|}{$-0.5$} &
\multicolumn{1}{|c|}{\color{gris}0.0904\color{black}} &
\multicolumn{1}{c|}{0.1973} & 0.2037 & \color{blue}0.0541\color{black} &
0.1252 & 0.1086\\
\multicolumn{1}{|c|}{} & \multicolumn{1}{|c|}{$0$} &
\multicolumn{1}{|c|}{\color{blue}0.0694\color{black}} &
\multicolumn{1}{c|}{\color{red}0.0761\color{black}} & \color{blue}%
0.0686\color{black} & \color{blue}0.0488\color{black} & \color{blue}%
0.0646\color{black} & \color{blue}0.0592\color{black}\\
\multicolumn{1}{|c|}{$0$} & \multicolumn{1}{|c|}{$2/3$} &
\multicolumn{1}{|c|}{\color{blue}0.0514\color{black}} &
\multicolumn{1}{c|}{\color{blue}0.0455\color{black}} & \color{blue}%
0.0408\color{black} & \color{blue}0.0443\color{black} & \color{blue}%
0.0431\color{black} & \color{blue}0.0450\color{black}\\
\multicolumn{1}{|c|}{} & \multicolumn{1}{|c|}{$1$} &
\multicolumn{1}{|c|}{\color{blue}0.0479\color{black}} &
\multicolumn{1}{c|}{\color{blue}0.0441\color{black}} & \color{blue}%
0.0401\color{black} & \color{blue}0.0432\color{black} & \color{blue}%
0.0391\color{black} & \color{blue}0.0422\color{black}\\
\multicolumn{1}{|c|}{} & \multicolumn{1}{|c|}{$2$} &
\multicolumn{1}{|c|}{\color{blue}0.0500\color{black}} &
\multicolumn{1}{c|}{\color{red}0.0707\color{black}} & \color{blue}%
0.0512\color{black} & \color{blue}0.0441\color{black} & \color{blue}%
0.0431\color{black} & \color{blue}0.0434\color{black}\\\hline
\multicolumn{1}{|c|}{} & \multicolumn{1}{|c|}{$-0.5$} &
\multicolumn{1}{|c|}{0.1407} & \multicolumn{1}{c|}{0.2595} & 0.2600 &
\color{blue}0.0642\color{black} & 0.1374 & 0.1292\\
\multicolumn{1}{|c|}{} & \multicolumn{1}{|c|}{$0$} &
\multicolumn{1}{|c|}{\color{red}0.0808\color{black}} &
\multicolumn{1}{c|}{\color{red}0.0807\color{black}} & \color{red}%
0.0705\color{black} & \color{blue}0.0523\color{black} & \color{blue}%
0.0670\color{black} & \color{blue}0.0603\color{black}\\
\multicolumn{1}{|c|}{$2/3$} & \multicolumn{1}{|c|}{$2/3$} &
\multicolumn{1}{|c|}{\color{blue}0.0376\color{black}} &
\multicolumn{1}{c|}{\color{red}0.0257\color{black}} & 0.0244 & \color
{blue}0.0411\color{black} & \color{blue}0.0372\color{black} & \color
{red}0.0350\color{black}\\
\multicolumn{1}{|c|}{} & \multicolumn{1}{|c|}{$1$} &
\multicolumn{1}{|c|}{\color{red}0.0254\color{black}} &
\multicolumn{1}{c|}{0.0188} & 0.0194 & \color{blue}0.0369\color{black} &
\color{red}0.0308\color{black} & \color{red}0.0295\color{black}\\
\multicolumn{1}{|c|}{} & \multicolumn{1}{|c|}{$2$} &
\multicolumn{1}{|c|}{0.0118} & \multicolumn{1}{c|}{0.0173} & 0.0159 &
\color{red}0.0280\color{black} & 0.0214 & 0.0243\\\hline
\multicolumn{1}{|c|}{} & \multicolumn{1}{|c|}{$-0.5$} &
\multicolumn{1}{|c|}{0.1725} & \multicolumn{1}{c|}{0.2872} & 0.2832 &
\color{red}0.0743\color{black} & 0.1386 & 0.1500\\
\multicolumn{1}{|c|}{} & \multicolumn{1}{|c|}{$0$} &
\multicolumn{1}{|c|}{\color{red}0.0923\color{black}} &
\multicolumn{1}{c|}{\color{red}0.0883\color{black}} & \color{red}%
0.0780\color{black} & \color{blue}0.0578\color{black} & \color{blue}%
0.0658\color{black} & \color{red}0.0726\color{black}\\
\multicolumn{1}{|c|}{$1$} & \multicolumn{1}{|c|}{$2/3$} &
\multicolumn{1}{|c|}{\color{blue}0.0378\color{black}} &
\multicolumn{1}{c|}{\color{red}0.0257\color{black}} & \color{red}%
0.0298\color{black} & \color{blue}0.0422\color{black} & \color{red}%
0.0347\color{black} & \color{blue}0.0371\color{black}\\
\multicolumn{1}{|c|}{} & \multicolumn{1}{|c|}{$1$} &
\multicolumn{1}{|c|}{0.0237} & \multicolumn{1}{c|}{0.0197} & \color
{red}0.0258\color{black} & \color{blue}0.0364\color{black} & \color
{red}0.0288\color{black} & \color{red}0.0290\color{black}\\
\multicolumn{1}{|c|}{} & \multicolumn{1}{|c|}{$2$} &
\multicolumn{1}{|c|}{0.0084} & \multicolumn{1}{c|}{0.0181} & \color
{red}0.0295\color{black} & 0.0239 & 0.0218 & 0.0179\\\hline
\end{tabular}
$
\end{center}
\end{table}\begin{table}[ptbh]
\caption{$\alpha_{n}^{(\lambda_{1},\lambda_{2})}$ with $n\in\{400,550\}$}%
\label{Table3}
\begin{center}
$%
\begin{tabular}
[c]{cc|ccc|ccc|}\cline{3-8}
&  & \multicolumn{3}{|c|}{$n=400$} & \multicolumn{3}{|c|}{$n=550$}\\\hline
\multicolumn{1}{|c|}{$\lambda_{2}$} & \multicolumn{1}{|c|}{$\lambda_{1}$} &
\multicolumn{1}{|c|}{(\ref{eq43})} & \multicolumn{1}{c|}{(\ref{eq30}%
)-(\ref{eq42})} & (\ref{eq44})-(\ref{eq45}) & (\ref{eq43}) &
\multicolumn{1}{|c}{(\ref{eq30})-(\ref{eq42})} &
\multicolumn{1}{|c|}{(\ref{eq44})-(\ref{eq45})}\\\hline
\multicolumn{1}{|c|}{} & \multicolumn{1}{|c|}{$-0.5$} &
\multicolumn{1}{|c|}{\color{blue}0.0587\color{black}} &
\multicolumn{1}{c|}{\color{red}0.0800\color{black}} &
\multicolumn{1}{c|}{\color{red}0.0704\color{black}} &
\multicolumn{1}{|c|}{\color{blue}0.0554\color{black}} &
\multicolumn{1}{c|}{\color{blue}0.0664\color{black}} & \color{blue}%
0.0664\color{black}\\
\multicolumn{1}{|c|}{} & \multicolumn{1}{|c|}{$0$} &
\multicolumn{1}{|c|}{\color{blue}0.0556\color{black}} &
\multicolumn{1}{c|}{\color{blue}0.0589\color{black}} &
\multicolumn{1}{c|}{\color{blue}0.0553\color{black}} &
\multicolumn{1}{|c|}{\color{blue}0.0530\color{black}} &
\multicolumn{1}{c|}{\color{blue}0.0566\color{black}} & \color{blue}%
0.0566\color{black}\\
\multicolumn{1}{|c|}{$0$} & \multicolumn{1}{|c|}{$2/3$} &
\multicolumn{1}{|c|}{\color{blue}0.0526\color{black}} &
\multicolumn{1}{c|}{\color{blue}0.0479\color{black}} &
\multicolumn{1}{c|}{\color{blue}0.0482\color{black}} &
\multicolumn{1}{|c|}{\color{blue}0.0515\color{black}} &
\multicolumn{1}{c|}{\color{blue}0.0495\color{black}} & \color{blue}%
0.0495\color{black}\\
\multicolumn{1}{|c|}{} & \multicolumn{1}{|c|}{$1$} &
\multicolumn{1}{|c|}{\color{blue}0.0518\color{black}} &
\multicolumn{1}{c|}{\color{blue}0.0454\color{black}} &
\multicolumn{1}{c|}{\color{blue}0.0461\color{black}} &
\multicolumn{1}{|c|}{\color{blue}0.0510\color{black}} &
\multicolumn{1}{c|}{\color{blue}0.0483\color{black}} & \color{blue}%
0.0483\color{black}\\
\multicolumn{1}{|c|}{} & \multicolumn{1}{|c|}{$2$} &
\multicolumn{1}{|c|}{\color{blue}0.0520\color{black}} &
\multicolumn{1}{c|}{\color{blue}0.0460\color{black}} &
\multicolumn{1}{c|}{\color{blue}0.0464\color{black}} &
\multicolumn{1}{|c|}{\color{blue}0.0509\color{black}} &
\multicolumn{1}{c|}{\color{blue}0.0483\color{black}} & \color{blue}%
0.0483\color{black}\\\hline
\multicolumn{1}{|c|}{} & \multicolumn{1}{|c|}{$-0.5$} &
\multicolumn{1}{|c|}{\color{blue}0.0585\color{black}} &
\multicolumn{1}{c|}{\color{red}0.0877\color{black}} &
\multicolumn{1}{c|}{\color{red}0.0800\color{black}} &
\multicolumn{1}{|c|}{\color{blue}0.0578\color{black}} &
\multicolumn{1}{c|}{\color{red}0.0703\color{black}} & \color{red}%
0.0703\color{black}\\
\multicolumn{1}{|c|}{} & \multicolumn{1}{|c|}{$0$} &
\multicolumn{1}{|c|}{\color{blue}0.0503\color{black}} &
\multicolumn{1}{c|}{\color{blue}0.0609\color{black}} &
\multicolumn{1}{c|}{\color{blue}0.0555\color{black}} &
\multicolumn{1}{|c|}{\color{blue}0.0541\color{black}} &
\multicolumn{1}{c|}{\color{blue}0.0585\color{black}} & \color{blue}%
0.0585\color{black}\\
\multicolumn{1}{|c|}{$2/3$} & \multicolumn{1}{|c|}{$2/3$} &
\multicolumn{1}{|c|}{\color{blue}0.0441\color{black}} &
\multicolumn{1}{c|}{\color{blue}0.0440\color{black}} &
\multicolumn{1}{c|}{\color{blue}0.0413\color{black}} &
\multicolumn{1}{|c|}{\color{blue}0.0494\color{black}} &
\multicolumn{1}{c|}{\color{blue}0.0463\color{black}} & \color{blue}%
0.0463\color{black}\\
\multicolumn{1}{|c|}{} & \multicolumn{1}{|c|}{$1$} &
\multicolumn{1}{|c|}{\color{blue}0.0415\color{black}} &
\multicolumn{1}{c|}{\color{blue}0.0384\color{black}} &
\multicolumn{1}{c|}{\color{blue}0.0372\color{black}} &
\multicolumn{1}{|c|}{\color{blue}0.0480\color{black}} &
\multicolumn{1}{c|}{\color{blue}0.0418\color{black}} & \color{blue}%
0.0418\color{black}\\
\multicolumn{1}{|c|}{} & \multicolumn{1}{|c|}{$2$} &
\multicolumn{1}{|c|}{\color{blue}0.0375\color{black}} &
\multicolumn{1}{c|}{\color{red}0.0299\color{black}} &
\multicolumn{1}{c|}{\color{red}0.0309\color{black}} &
\multicolumn{1}{|c|}{\color{blue}0.0433\color{black}} &
\multicolumn{1}{c|}{\color{red}0.0353\color{black}} & \color{red}%
0.0353\color{black}\\\hline
\multicolumn{1}{|c|}{} & \multicolumn{1}{|c|}{$-0.5$} &
\multicolumn{1}{|c|}{\color{blue}0.0628\color{black}} &
\multicolumn{1}{c|}{\color{red}0.0887\color{black}} &
\multicolumn{1}{c|}{0.0960} & \multicolumn{1}{|c|}{\color{blue}0.0609\color
{black}} & \multicolumn{1}{c|}{\color{red}0.0738\color{black}} & \color
{red}0.0738\color{black}\\
\multicolumn{1}{|c|}{} & \multicolumn{1}{|c|}{$0$} &
\multicolumn{1}{|c|}{\color{blue}0.0541\color{black}} &
\multicolumn{1}{c|}{\color{blue}0.0581\color{black}} &
\multicolumn{1}{c|}{\color{blue}0.0653\color{black}} &
\multicolumn{1}{|c|}{\color{blue}0.0559\color{black}} &
\multicolumn{1}{c|}{\color{blue}0.0584\color{black}} & \color{blue}%
0.0584\color{black}\\
\multicolumn{1}{|c|}{$1$} & \multicolumn{1}{|c|}{$2/3$} &
\multicolumn{1}{|c|}{\color{blue}0.0444\color{black}} &
\multicolumn{1}{c|}{\color{blue}0.0413\color{black}} &
\multicolumn{1}{c|}{\color{blue}0.0439\color{black}} &
\multicolumn{1}{|c|}{\color{blue}0.0500\color{black}} &
\multicolumn{1}{c|}{\color{blue}0.0457\color{black}} & \color{blue}%
0.0457\color{black}\\
\multicolumn{1}{|c|}{} & \multicolumn{1}{|c|}{$1$} &
\multicolumn{1}{|c|}{\color{blue}0.0411\color{black}} &
\multicolumn{1}{c|}{\color{blue}0.0367\color{black}} &
\multicolumn{1}{c|}{\color{blue}0.0371\color{black}} &
\multicolumn{1}{|c|}{\color{blue}0.0479\color{black}} &
\multicolumn{1}{c|}{\color{blue}0.0415\color{black}} & \color{blue}%
0.0415\color{black}\\
\multicolumn{1}{|c|}{} & \multicolumn{1}{|c|}{$2$} &
\multicolumn{1}{|c|}{\color{red}0.0344\color{black}} &
\multicolumn{1}{c|}{\color{red}0.0290\color{black}} &
\multicolumn{1}{c|}{\color{red}0.0256\color{black}} &
\multicolumn{1}{|c|}{\color{blue}0.0410\color{black}} &
\multicolumn{1}{c|}{\color{red}0.0335\color{black}} & \color{red}%
0.0335\color{black}\\\hline
\end{tabular}
$
\end{center}
\end{table}To illustrate some representative values of the powers, in Table
\ref{Table4} the simulated exact powers focussed on the tests (\ref{eq44}%
)-(\ref{eq45}) are shown for the considered $12$ points. The so-called
\emph{size corrected average gradient},\ defined as%
\[
\gamma_{n}^{(\lambda_{1},\lambda_{2})}=\tfrac{\left(  \frac{1}{12}\left(
{\textstyle\sum\limits_{i=1}^{6}}
\left(  \frac{\beta_{n}^{(\lambda_{1},\lambda_{2})}(i)-\alpha_{n}%
^{(\lambda_{1},\lambda_{2})}}{\delta_{1}(i)+\delta_{2}(i)}\right)  ^{\!2}+%
{\textstyle\sum\limits_{i=7}^{12}}
\left(  \frac{\beta_{n}^{(\lambda_{1},\lambda_{2})}(i)-\alpha_{n}%
^{(\lambda_{1},\lambda_{2})}}{\beta(i)}\right)  ^{2}\right)  \right)
^{\frac{1}{2}}}{\alpha_{n}^{(\lambda_{1},\lambda_{2})}},
\]
is an overall measure of performance of the simulated exact size as well as
the simulated exact powers (this measure was introduced for the first time in
Rivas et al. \cite{Rivas}). Such a measure is interpreted a normalized mean
rate of power gain with respect to the null hypothesis along the considered
alternatives and it is therefore useful as criterion to select a test
statistic ($\lambda_{1}$) as well as its estimator ($\lambda_{2}$) with the
maximum value of $\gamma_{n}^{(\lambda_{1},\lambda_{2})}$. The values of
$\gamma_{n}^{(\lambda_{1},\lambda_{2})}$\ for the same kind of test-statistics
considered in Tables \ref{Table2} and \ref{Table3}\ are shown in Tables
\ref{Table5} and \ref{Table6}. \begin{table}[ptbh]
\caption{$\beta_{n}^{(\lambda_{1},\lambda_{2})}(i) $ with $\lambda_{2}%
=\frac{2}{3}$ for tests (\ref{eq44})-(\ref{eq45})}%
\label{Table4}
\begin{center}
$%
\begin{tabular}
[c]{cccccccc}\hline
$n$ & $\lambda_{1}$ & $i=1$ & $i=4$ & $i=2$ & $i=5$ & $i=3$ & $i=6$\\\hline
& $0$ & 0.0720 & 0.0749 & 0.1139 & 0.1201 & 0.1749 & 0.1745\\
$100$ & $2/3$ & 0.0188 & 0.0163 & 0.0382 & 0.0389 & 0.0736 & 0.0745\\
& $1$ & 0.0110 & 0.0088 & 0.0251 & 0.0252 & 0.0522 & 0.0502\\
& $2$ & 0.0038 & 0.0035 & 0.0095 & 0.0111 & 0.0271 & 0.0257\\\hline
& $0$ & 0.0996 & 0.0989 & 0.2123 & 0.2091 & 0.3669 & 0.3637\\
$250$ & $2/3$ & 0.0582 & 0.0571 & 0.1522 & 0.1516 & 0.2956 & 0.2945\\
& $1$ & 0.0483 & 0.0468 & 0.1320 & 0.1330 & 0.2737 & 0.2716\\
& $2$ & 0.0337 & 0.0323 & 0.1063 & 0.1047 & 0.2303 & 0.2329\\\hline
& $0$ & 0.1251 & 0.1161 & 0.3162 & 0.3230 & 0.5575 & 0.5631\\
$400$ & $2/3$ & 0.0973 & 0.0882 & 0.2760 & 0.2837 & 0.5226 & 0.5249\\
& $1$ & 0.0859 & 0.0797 & 0.2616 & 0.2692 & 0.5072 & 0.5093\\
& $2$ & 0.0698 & 0.0639 & 0.2334 & 0.2427 & 0.4786 & 0.4840\\\hline
& $0$ & 0.1536 & 0.1453 & 0.4363 & 0.4410 & 0.7292 & 0.7347\\
$550$ & $2/3$ & 0.1337 & 0.1202 & 0.4072 & 0.4102 & 0.7072 & 0.7178\\
& $1$ & 0.1258 & 0.1096 & 0.3942 & 0.3984 & 0.6995 & 0.7074\\
& $2$ & 0.1095 & 0.0941 & 0.3734 & 0.3757 & 0.6832 & 0.6911\\\hline
$n$ & $\lambda_{1}$ & $i=7$ & $i=10$ & $i=8$ & $i=11$ & $i=9$ & $i=12$\\\hline
& $0$ & 0.0678 & 0.0724 & 0.1197 & 0.1370 & 0.2419 & 0.2880\\
$100$ & $2/3$ & 0.0572 & 0.0634 & 0.1050 & 0.1206 & 0.2183 & 0.2646\\
& $1$ & 0.0527 & 0.0597 & 0.0993 & 0.1148 & 0.2088 & 0.2562\\
& $2$ & 0.0511 & 0.0550 & 0.0931 & 0.1061 & 0.1949 & 0.2452\\\hline
& $0$ & 0.1798 & 0.1980 & 0.3597 & 0.3901 & 0.6588 & 0.7164\\
$250$ & $2/3$ & 0.1704 & 0.1861 & 0.3477 & 0.3739 & 0.6432 & 0.7010\\
& $1$ & 0.1661 & 0.1812 & 0.3418 & 0.3677 & 0.6363 & 0.6955\\
& $2$ & 0.1592 & 0.1750 & 0.3311 & 0.3557 & 0.6252 & 0.6859\\\hline
& $0$ & 0.3112 & 0.3369 & 0.5785 & 0.6124 & 0.8739 & 0.9190\\
$400$ & $2/3$ & 0.3037 & 0.3258 & 0.5683 & 0.6041 & 0.8685 & 0.9153\\
& $1$ & 0.3000 & 0.3229 & 0.5636 & 0.6009 & 0.8664 & 0.9137\\
& $2$ & 0.2917 & 0.3156 & 0.5544 & 0.5935 & 0.8610 & 0.9098\\\hline
& $0$ & 0.4249 & 0.4578 & 0.7272 & 0.7734 & 0.9630 & 0.9806\\
$550$ & $2/3$ & 0.4177 & 0.4510 & 0.7215 & 0.7671 & 0.9607 & 0.9798\\
& $1$ & 0.4147 & 0.4480 & 0.7191 & 0.7651 & 0.9598 & 0.9797\\
& $2$ & 0.4077 & 0.4431 & 0.7117 & 0.7593 & 0.9581 & 0.9789\\\hline
\end{tabular}
$
\end{center}
\end{table}\begin{table}[ptbh]
\caption{$\gamma_{n}^{(\lambda_{1},\lambda_{2})}$ with $n\in\{100,250\}$}%
\label{Table5}
\begin{center}
$%
\begin{tabular}
[c]{ccc|cc|ccc|}\cline{3-8}
&  & \multicolumn{3}{|c|}{$n=100$} & \multicolumn{3}{|c|}{$n=250$}\\\hline
\multicolumn{1}{|c}{$\lambda_{2}$} & \multicolumn{1}{|c}{$\lambda_{1}$} &
\multicolumn{1}{|c|}{(\ref{eq43})} & \multicolumn{1}{|c|}{(\ref{eq30}%
)-(\ref{eq42})} & (\ref{eq44})-(\ref{eq45}) & (\ref{eq43}) &
\multicolumn{1}{|c|}{(\ref{eq30})-(\ref{eq42})} & (\ref{eq44})-(\ref{eq45}%
)\\\hline
\multicolumn{1}{|c|}{} & \multicolumn{1}{c|}{$0$} & 1.9498 & 0.8109 &
\multicolumn{1}{|c|}{1.8522} & \multicolumn{1}{|c|}{8.8777} &
\multicolumn{1}{c|}{2.2680} & 6.5096\\
\multicolumn{1}{|c|}{$0$} & \multicolumn{1}{c|}{$2/3$} & 2.3191 & 0.8173 &
\multicolumn{1}{|c|}{3.3489} & \multicolumn{1}{|c|}{9.5346} &
\multicolumn{1}{c|}{3.4695} & 8.5937\\
\multicolumn{1}{|c|}{} & \multicolumn{1}{c|}{$1$} & 2.3989 & 0.8030 &
\multicolumn{1}{|c|}{3.3222} & \multicolumn{1}{|c|}{9.7044} &
\multicolumn{1}{c|}{3.8527} & 9.2039\\
\multicolumn{1}{|c|}{} & \multicolumn{1}{c|}{$2$} & 2.3323 & 0.7436 &
\multicolumn{1}{|c|}{2.5181} & \multicolumn{1}{|c|}{9.4949} &
\multicolumn{1}{c|}{3.5113} & 9.0488\\\hline
\multicolumn{1}{|c|}{} & \multicolumn{1}{c|}{$0$} & 2.2460 & 0.2623 &
\multicolumn{1}{|c|}{1.4553} & 8.4591 & \multicolumn{1}{|c}{2.0100} &
\multicolumn{1}{|c|}{6.5121}\\
\multicolumn{1}{|c|}{$2/3$} & \multicolumn{1}{c|}{$2/3$} & 3.3642 & 1.2907 &
\multicolumn{1}{|c|}{4.5181} & 10.0525 & \multicolumn{1}{|c}{3.4108} &
\multicolumn{1}{|c|}{11.1638}\\
\multicolumn{1}{|c|}{} & \multicolumn{1}{c|}{$1$} & 4.1761 & 1.8014 &
\multicolumn{1}{|c|}{5.5197} & 10.8193 & \multicolumn{1}{|c}{3.9618} &
\multicolumn{1}{|c|}{13.1137}\\
\multicolumn{1}{|c|}{} & \multicolumn{1}{c|}{$2$} & 6.1807 & 1.6047 &
\multicolumn{1}{|c|}{6.4614} & 13.1578 & \multicolumn{1}{|c}{5.3177} &
\multicolumn{1}{|c|}{15.5109}\\\hline
\multicolumn{1}{|c|}{} & \multicolumn{1}{c|}{$0$} & 2.0916 & 0.8894 &
\multicolumn{1}{|c|}{1.2676} & 7.8486 & \multicolumn{1}{|c}{1.7541} &
\multicolumn{1}{|c|}{5.5833}\\
\multicolumn{1}{|c|}{$1$} & \multicolumn{1}{c|}{$2/3$} & 3.3562 & 0.8595 &
\multicolumn{1}{|c|}{3.3604} & 9.8390 & \multicolumn{1}{|c}{3.3627} &
\multicolumn{1}{|c|}{10.5545}\\
\multicolumn{1}{|c|}{} & \multicolumn{1}{c|}{$1$} & 4.2558 & 0.9214 &
\multicolumn{1}{|c|}{3.7848} & 10.8838 & \multicolumn{1}{|c}{4.1772} &
\multicolumn{1}{|c|}{12.5275}\\
\multicolumn{1}{|c|}{} & \multicolumn{1}{c|}{$2$} & 7.1583 & 1.0689 &
\multicolumn{1}{|c|}{3.0949} & 14.5192 & \multicolumn{1}{|c}{6.2040} &
\multicolumn{1}{|c|}{15.9077}\\\hline
\end{tabular}
$
\end{center}
\end{table}\begin{table}[ptbh]
\caption{$\gamma_{n}^{(\lambda_{1},\lambda_{2})}$ with $n\in\{400,550\}$}%
\label{Table6}
\begin{center}
$%
\begin{tabular}
[c]{ccccc|ccc|}\cline{3-8}
&  & \multicolumn{3}{|c|}{$n=400$} & \multicolumn{3}{|c|}{$n=550$}\\\hline
\multicolumn{1}{|c}{$\lambda_{2}$} & \multicolumn{1}{|c}{$\lambda_{1}$} &
\multicolumn{1}{|c}{(\ref{eq43})} & \multicolumn{1}{|c|}{(\ref{eq30}%
)-(\ref{eq42})} & (\ref{eq44})-(\ref{eq45}) & \multicolumn{1}{|c|}{(\ref{eq43}%
)} & \multicolumn{1}{|c|}{(\ref{eq30})-(\ref{eq42})} & (\ref{eq44}%
)-(\ref{eq45})\\\hline
\multicolumn{1}{|c|}{} & \multicolumn{1}{c|}{$0$} &
\multicolumn{1}{c|}{11.5475} & \multicolumn{1}{c|}{5.1013} & 10.5516 &
\multicolumn{1}{|c|}{15.1129} & \multicolumn{1}{c|}{7.5141} & 13.4170\\
\multicolumn{1}{|c|}{$0$} & \multicolumn{1}{c|}{$2/3$} &
\multicolumn{1}{c|}{12.0922} & \multicolumn{1}{c|}{6.3045} & 12.1017 &
\multicolumn{1}{|c|}{15.4706} & \multicolumn{1}{c|}{8.6282} & 14.9108\\
\multicolumn{1}{|c|}{} & \multicolumn{1}{c|}{$1$} &
\multicolumn{1}{c|}{12.2447} & \multicolumn{1}{c|}{6.6720} & 12.6871 &
\multicolumn{1}{|c|}{15.5996} & \multicolumn{1}{c|}{8.8334} & 15.0571\\
\multicolumn{1}{|c|}{} & \multicolumn{1}{c|}{$2$} &
\multicolumn{1}{c|}{12.1903} & \multicolumn{1}{c|}{6.6276} & 12.6977 &
\multicolumn{1}{|c|}{15.6361} & \multicolumn{1}{c|}{8.8885} & 14.9610\\\hline
\multicolumn{1}{|c|}{} & \multicolumn{1}{c|}{$0$} & 13.0615 &
\multicolumn{1}{|c}{4.8961} & \multicolumn{1}{|c|}{10.9942} & 14.8684 &
\multicolumn{1}{|c}{7.3124} & \multicolumn{1}{|c|}{13.7563}\\
\multicolumn{1}{|c|}{$2/3$} & \multicolumn{1}{c|}{$2/3$} & 14.4991 &
\multicolumn{1}{|c}{6.3327} & \multicolumn{1}{|c|}{14.7299} & 16.0838 &
\multicolumn{1}{|c}{8.8237} & \multicolumn{1}{|c|}{17.1009}\\
\multicolumn{1}{|c|}{} & \multicolumn{1}{c|}{$1$} & 15.2220 &
\multicolumn{1}{|c}{7.1010} & \multicolumn{1}{|c|}{16.2884} & 16.4405 &
\multicolumn{1}{|c}{9.5901} & \multicolumn{1}{|c|}{18.3816}\\
\multicolumn{1}{|c|}{} & \multicolumn{1}{c|}{$2$} & 16.3217 &
\multicolumn{1}{|c}{8.7403} & \multicolumn{1}{|c|}{19.3697} & 18.0044 &
\multicolumn{1}{|c}{10.9238} & \multicolumn{1}{|c|}{21.2854}\\\hline
\multicolumn{1}{|c|}{} & \multicolumn{1}{c|}{$0$} & 12.2295 &
\multicolumn{1}{|c}{4.3335} & \multicolumn{1}{|c|}{10.1186} & 14.4705 &
\multicolumn{1}{|c}{6.8205} & \multicolumn{1}{|c|}{12.7432}\\
\multicolumn{1}{|c|}{$1$} & \multicolumn{1}{c|}{$2/3$} & 14.4409 &
\multicolumn{1}{|c}{6.3769} & \multicolumn{1}{|c|}{14.1149} & 15.8989 &
\multicolumn{1}{|c}{8.8589} & \multicolumn{1}{|c|}{16.2596}\\
\multicolumn{1}{|c|}{} & \multicolumn{1}{c|}{$1$} & 15.3240 &
\multicolumn{1}{|c}{7.4787} & \multicolumn{1}{|c|}{15.7674} & 16.4489 &
\multicolumn{1}{|c}{10.0339} & \multicolumn{1}{|c|}{17.8500}\\
\multicolumn{1}{|c|}{} & \multicolumn{1}{c|}{$2$} & 17.4356 &
\multicolumn{1}{|c}{10.5383} & \multicolumn{1}{|c|}{19.5188} & 18.8438 &
\multicolumn{1}{|c}{12.4853} & \multicolumn{1}{|c|}{21.8208}\\\hline
\end{tabular}
$
\end{center}
\end{table}

\subsection{Conclusions}

A clear conclusion from the results of Tables \ref{Table5} and \ref{Table6} is
the best performance of the sequence of tests (\ref{eq44})-(\ref{eq45})
compared with (\ref{eq30})-(\ref{eq42}), essentially because the tests
(\ref{eq44})-(\ref{eq45}) were much more powerful than (\ref{eq30}%
)-(\ref{eq42}). One possible reason, as it is explained in Agresti \cite[page
373]{Agresti}, could be that the power of a chi-square test tends to increase
when degrees of freedom decrease. This explanation is valid for (\ref{eq45}),
however even a greater value of degrees of freedom of (\ref{eq44}), the
computed powers of the test (\ref{eq44}) in this experiment were also greater
than the powers of the test (\ref{eq30}). In spite of that, we think that the
last result could be affected by the common choice of (\ref{eq46}), which
means that the points within the region $M_{\mathcal{QS}}-M_{\mathcal{OQS}}%
$\ are excluded in the alternative hypothesis of (\ref{eq44}) because these
points fall within the points which are included in the null hypothesis of
(\ref{eq30}). On the other hand if one desires to compare the test-statistics
associated with (\ref{eq44})-(\ref{eq45}) with respect to ones of
(\ref{eq43}), even perhaps better performance of (\ref{eq44})-(\ref{eq45}%
)\ there is no a big difference, and thus unless there is an evidence for
thinking that there exists OQS before carrying out a MH test, it would be more
convenient to use the unconditional HM test (\ref{eq43}). Within the values of
$\gamma_{n}^{(\lambda_{1},\lambda_{2})}$\ for the same test, either
(\ref{eq44})-(\ref{eq45})\ or (\ref{eq43}), the variability is greater for
(\ref{eq44})-(\ref{eq45}) because the simulated exact sizes are also more
variable. Apart from the criterion of the corrected average gradient, a
criterion for excluding from the study the simulated exact sizes which are not
close or fairly close to the nominal size ($\alpha=0.05$) is necessary.
Through the criterion given by Dale \cite{Dale},\ the inequality%
\[
\left\vert \operatorname{logit}(1-\alpha_{n}^{(\lambda_{1},\lambda_{2}%
)})-\operatorname{logit}(1-\alpha)\right\vert \leq\epsilon,
\]
with $\operatorname{logit}(p)\equiv\log(p/(1-p))$,\ is considered, so that the
two probabilities, $\alpha_{n}^{(\lambda_{1},\lambda_{2})}$ and $\alpha$, are
considered to be close\ if they satisfy such a inequality with $\epsilon
=0.35$\ and fairly close\ if they satisfy with $\epsilon=0.7$. Note that for
$\alpha=0.05$, $\epsilon=0.35$\ corresponds to $\alpha_{n}^{(\lambda
_{1},\lambda_{2})}\in\lbrack0.0357,0.0695]$ and $\epsilon=0.7$\ corresponds to
$\alpha_{n}^{(\lambda_{1},\lambda_{2})}\in\lbrack0.0254,0.0357)\cup
(0.0695,0.0954]$. Those simulated exact sizes which are taken as
close\ according to the Dale's criterion\ have been marked in blue color in
Tables \ref{Table1} and \ref{Table2}, and in red color fairly close\ simulated
exact sizes. Finally, it is concluded that the best overall choice for the
test statistics $T^{\phi_{(\lambda_{1})}}(\widehat{\boldsymbol{\theta}}%
^{\phi_{(\lambda_{2})}})$\ is $(\lambda_{1},\lambda_{2})\in\{(1,1),(1,\frac
{2}{3})\}$, however for the smallest sample size ($n=100$) the test statistic
associated with $\lambda_{1}=1$ is not a good choice and it is better
$(\lambda_{1},\lambda_{2})\in\{(\frac{2}{3},1),(\frac{2}{3},\frac{2}{3})\}$.

To finalize, we would like to comment that LMLC can have been dealt in this
paper in a more general setting by following generalized log-linear models,
$C\log\left(  A\boldsymbol{m}(\boldsymbol{\theta}\mathbf{)}\right)
\mathbf{=}\boldsymbol{X}\boldsymbol{\theta}$ (see Lang \cite{Lang0}\ and
references therein). With these models it would be possible to consider
loglinear constraints for the marginal distributions by considering
$C=\boldsymbol{I}_{k}$ and $A\neq\boldsymbol{I}_{k}$. Furthermore, using
minimum power divergence estimators a different type of application for
$\log\left(  A\boldsymbol{m}(\boldsymbol{\theta}\mathbf{)}\right)
\mathbf{=}\boldsymbol{X}\boldsymbol{\theta}$, with $A\neq\boldsymbol{I}_{k}$
and Poisson sampling, can be found in Mart\'{\i}n and Li \cite{MartinLi}.

\end{document}